\documentclass[11pt]{amsart}
\usepackage[mathscr]{euscript}
\usepackage{amssymb}

\usepackage{amscd}
\usepackage{amsmath, amssymb}
\usepackage{amsfonts}
\usepackage[colorlinks,linkcolor=blue,citecolor=blue, pdfstartview=FitH]{hyperref}
\usepackage[all]{xy}

\numberwithin{equation}{section}

\theoremstyle{plain}
\newtheorem{thm}{Theorem}[section]
\newtheorem{theorem}[thm]{Theorem}
\newtheorem{lemma}[thm]{Lemma}
\newtheorem{corollary}[thm]{Corollary}
\newtheorem{proposition}[thm]{Proposition}

\theoremstyle{definition}

\newtheorem{remark}[thm]{Remark}

\newtheorem{definition}[thm]{Definition}

\newtheorem{example}[thm]{Example}

\newtheorem{defn-thm}[thm]{Definition-Theorem}

\newcommand{\im}{{ \mathrm{im}\,}}

\newcommand{\C}{{ \mathbb{C} }}

\newcommand{\N}{{ \mathbb{N} }}

\newcommand{\p}{{ \partial }}
\newcommand{\pb}{{ \bar{\partial} }}

\newcommand{\wg}{{\,\wedge\,}}

\begin{document}
\title{Deformed Aeppli cohomology: canonical deformations and jumping formulas}
\author{Yan Hu, Wei Xia}
\address{Yan Hu, Mathematical Science Research Center, Chongqing University of Technology, Chongqing 400054, China} \email{huyan1128@126.com}
\address{Wei Xia, Mathematical Science Research Center, Chongqing University of Technology, Chongqing 400054, China} \email{xiaweiwei3@126.com, xiawei@cqut.edu.cn}

\thanks{This work is supported by the National Natural Science Foundation of China No. 11901590, the Natural Science Foundation of Chongqing (China) No. CSTB2022NSCQ-MSX0876. This work is also partially supported by the Natural Science Foundation of Chongqing (China) No. CSTB2024NSCQ-LZX0040 and No. CSTB2023NSCQ-LZX0042.}
\date{\today}

\begin{abstract}
Given a complex analytic family of complex manifolds, we consider canonical Aeppli deformations of $(p,q)$-forms and study its relations to the varying of dimension of the deformed Aeppli cohomology $\dim  H^{\bullet,\bullet}_{A\phi(t)}(X)$. In particular, we prove the jumping formula for the deformed Aeppli cohomology $H^{\bullet,\bullet}_{A\phi(t)}(X)$. As a direct consequence, $\dim  H^{p,q}_{A\phi(t)}(X)$ remains constant iff the Bott-Chern deformations of $(n-p,n-q)$-forms and the Aeppli deformations of $(n-p-1,n-q-1)$-forms are canonically unobstructed. Furthermore, the Bott-Chern/Aeppli deformations are shown to be unobstructed if some weak forms of $\p\pb$-lemma is satisfied.
\end{abstract}

\maketitle

\tableofcontents
\section{Introduction}
Let $\pi: (\mathcal{X}, X)\to (B,0)$ be a complex analytic family such that for each $t\in B$ the complex structure on $X_t:=\pi^{-1}(t)$ is represented by Beltrami differential $\phi(t)\in A^{0,1}(X,T_X^{1,0})$. Here and throughout this paper, we will always assume $B$ is a polydisc with sufficiently small radius and the reference fiber $X_0$ is denoted by $X$ unless otherwise stated. The deformed Bott-Chern cohomology and the deformed Aeppli cohomology is defined \cite{Xia19dBC} as follows
\begin{align}\label{eq-Aphit}
H^{p,q}_{BC\phi(t)}(X) :=& \frac{\ker d_{\phi(t)} \cap A^{p,q}(X)}{\im \p\pb_{\phi(t)} \cap A^{p,q}(X)},\\
 H^{p,q}_{A\phi(t)}(X) :=& \frac{\ker \p\pb_{\phi(t)} \cap A^{p,q}(X)}{(\im \p+ \im \pb_{\phi(t)}) \cap A^{p,q}(X)}.
\end{align}
Note that when $t=0$ these are just the usual Bott-Chern and Aeppli cohomology, see \cite{Ang13,AT13,AT15b,ADT16,RZ24,Sch07,YY17} and \cite{Pop15,PT20} for some recent works. One of the main motivations to study these deformed cohomology groups is that they contain information about the cohomology groups on the deformed manifold $X_t$ and thus could be used as a possible tool to prove deformation invariance of the dimension of corresponding cohomology groups (see \cite{LRY15,LSY09,RZ15,RZ18,RZ22,WZ20,Xia19dDol} for related discussions). Indeed, for any $p\in \mathbb{N}$
\[
H^{p,0}_{BC\phi(t)}(X)\cong H^{p,0}_{BC}(X_t),~H^{p,n}_{A\phi(t)}(X)\cong H^{p,n}_{A}(X_t), \quad \text{for~any}~t\in B.
\]
While the deformed Bott-Chern cohomology is already studied in \cite{Xia19dBC} where a deformation theory of Bott-Chern cohomology has been established, the deformed Aeppli cohomology is relatively unexplored. This is partly due to the fact that the deformed Aeppli cohomology is dual, via the Hodge star operator $*$, to the deformed Bott-Chern cohomology (see \eqref{eq-BC-A-dual}):
\[
*:H^{p,q}_{A\phi(t)}(X)\longrightarrow H^{n-p,n-q}_{BC\phi(t)}(X), \quad \text{for~any}~t\in B.
\]
It then seems reasonable to expect statements (e.g. the corresponding deformation theory) that is valid for the deformed Bott-Chern cohomology are somehow dual, and thus can be easily translated, to the case of deformed Aeppli cohomology. However, our results in this paper will show this is not entirely correct at least when the deformation theory is concerned.

In this paper, we will establish a deformation theory for the Aeppli cohomology. In particular, the relations between jumping of $\dim H^{p,q}_{A\phi(t)}(X)$ and canonical Aeppli deformations of $(p,q)$-forms will be studied.
\begin{definition}\label{def-deformation-Dol/BC/A-0}
Given a $(p,q)$-form
\[
y\in \ker\pb/\ker d/\ker\p\pb
\]
and $T\subseteq B$, which is an analytic subset of $B$ containing $0$,
\[
\text{a~\emph{Dolbeault/Bott-Chern/Aeppli deformation}~of~$y$~w.r.t.}~\pi: (\mathcal{X}, X)\to (B,0)~\text{on}~T
\]
is a family of $(p,q)$-forms $\sigma (t)$ such that
\begin{itemize}
  \item[1.] $\sigma (t)$ is holomorphic in $t\in T$ and $\sigma (0) = y$;
  \item[2.] For any $t\in T$, the following holds respectively:
  \begin{itemize} \item $\pb_{\phi(t)}\sigma (t) = 0$ (Dolbeault);
  \item $d_{\phi(t)}\sigma (t) = 0$ (Bott-Chern);
  \item $\p\pb_{\phi(t)}\sigma (t) = 0$ (Aeppli).
  \end{itemize}
\end{itemize}
\end{definition}
It turns out in order to capture the jumping of $\dim H^{p,q}_{A\phi(t)}(X)$ both the Bott-Chern deformations of $(n-p,n-q)$-forms and the Aeppli deformations of $(n-p-1,n-q-1)$-forms are involved. Our new jumping formulas provide a precise description about this dependence relations. More precisely, we have
\begin{theorem}[=Theorem \ref{thm-BCandA-jump-formula}]\label{thm-BCandA-jump-formula-0}
Let $\pi: (\mathcal{X}, X)\to (B,0)$ be a complex analytic family such that for each $t\in B$ the complex structure on $X_t$ is represented by Beltrami differential $\phi(t)$. Suppose $X$ is equipped with a fixed Hermitian metric. For each $(p,q)\in \mathbb{N}\times \mathbb{N}$, set\footnote{The convention is that $v^{p,q}_t=w^{p,q}_t=0$ whenever $p<0$ or $q<0$.}
\[
v^{p,q}_t:=\dim H_{BC}^{p,q}(X)-\dim \ker d_{\phi(t)}\cap\ker(\p\bar{\p})^*\cap A^{p,q}(X)\geq 0,
\]
and
\[
w^{p,q}_t:=\dim H_{A}^{p,q}(X)-\dim\ker\p\pb_{\phi(t)}\cap\ker d^*\cap A^{p,q}(X)\geq 0,
\]
then we have the jumping formula for the deformed Bott-Chern cohomology:
\begin{equation}\label{eq-BCjump-0}
\dim H_{BC}^{p,q}(X)=\dim H_{BC\phi(t)}^{p,q}(X)+v^{p,q}_t+w^{p-1,q-1}_t,
\end{equation}
and the jumping formula for the deformed Aeppli cohomology:
\begin{equation}\label{eq-Ajump-0}
\dim H_{A}^{p,q}(X)=\dim H_{A\phi(t)}^{p,q}(X)+v^{n-p,n-q}_t+w^{n-p-1,n-q-1}_t.
\end{equation}
\end{theorem}
In a previous work by the second author, a formula similar to \eqref{eq-BCjump-0} was proved \cite[Thm.\,4.13]{Xia19dBC} where $w^{p,q}_t$ is replaced by
\[
u^{p,q}_t=\dim H_{BC}^{p,q}(X)-\dim \ker \p\bar{\p}_{\phi(t)}\cap\left(\mathcal{H}_{BC}^{p,q}(X)+\im(\p\bar{\p})^*\right)\cap A^{p,q}(X).
\]
A major advantage of the new jumping formulas in Theorem \ref{thm-BCandA-jump-formula-0} lies in the fact that in contrast to $u^{p,q}_t$, both $v^{p,q}_t$ and $w^{p,q}_t$ have precise geometric meanings. Roughly speaking, $w^{p,q}_t$ ($v^{p,q}_t$) measures the dimension of the space of Aeppli (Bott-Chern) cohomology classes whose canonical deformations fail to be unobstructed at $t$, see \eqref{eq-wt-Wt}. In particular, Bott-Chern (Aeppli) deformations of $(p,q)$-forms are unobstructed iff $v^{p,q}_t\equiv0$ ($w^{p,q}_t\equiv0$). Because of this, we have the following
\begin{corollary}\label{coro-0-1}
In the situation of Theorem \ref{thm-BCandA-jump-formula-0}.
\begin{enumerate}
  \item The deformed Bott-Chern number $h_{BC\phi(t)}^{p,q}$ is independent of $t\in B$ if and only if the Bott-Chern deformations of $(p,q)$-forms and the Aeppli deformations of $(p-1,q-1)$-forms are canonically unobstructed;
  \item The deformed Aeppli number $h_{A\phi(t)}^{p,q}$ is independent of $t\in B$ if and only if the Bott-Chern deformations of $(n-p,n-q)$-forms and the Aeppli deformations of $(n-p-1,n-q-1)$-forms are canonically unobstructed;
\end{enumerate}
\end{corollary}
This gives a criterion for the deformation invariance of the deformed Bott-Chern and Aeppli numbers in terms of Bott-Chern/Aeppli deformations of forms. This picture is quite similar to results in \cite{LY15} where Lin-Ye discovered that the obstruction to extend Bott-Chern/Aeppli cohomology classes is not zero if the Bott-Chern/Aeppli numbers jumps. By using Corollary \ref{coro-0-1}, we are able to show the following (see \cite[Sec.\,5]{RWZ19} for related discussions):
\begin{theorem}[=Theorem \ref{thm-unobstructed-deformations}]
Let $\pi: (\mathcal{X}, X)\to (B,0)$ be a complex analytic family with Beltrami differentials $\phi(t)$. Suppose $X$ is equipped with a fixed Hermitian metric.
\begin{enumerate}
  \item Assume $\p_{A,\pb(\ker\p)}^{p-1,q+1}=0$. Then the Bott-Chern deformations of $(p,q)$-forms are canonically unobstructed;
  \item Assume $\p_{A,BC}^{p,q+1}=0$. Then the Aeppli deformations of $(p,q)$-forms are canonically unobstructed;
  \item Assume $\p_{A,\pb(\ker\p)}^{p-1,q+1}=0$ and $\p_{A,BC}^{p-1,q}=0$. Then the deformed Bott-Chern number $h_{BC\phi(t)}^{p,q}$ is independent of $t\in B$;
  \item Assume $\p_{A,\pb(\ker\p)}^{n-p-1,n-q+1}=0$ and $\p_{A,BC}^{n-p-1,n-q}=0$. Then the deformed Aeppli number $h_{A\phi(t)}^{p,q}$ is independent of $t\in B$.
\end{enumerate}
\end{theorem}
\begin{corollary}\label{coro-0}
The Bott-Chern number $\dim H_{BC}^{p,0}(X_t)=\dim H_{A}^{n,n-p}(X_t)$ is independent of $t\in B$ if $\p_{A,\pb(\ker\p)}^{p-1,1}=0$.
\end{corollary}
The conditions $\p_{A,\pb(\ker\p)}^{p-1,q+1}=0$ and $\p_{A,BC}^{p,q+1}=0$ are related to weak forms of $\p\pb$-lemma whose variants appear in many different context (see e.g. \cite{RWZ21,Pop19,Ale19,RZ18,AU17,AU16,FY11}. More precisely, $\p_{A,\pb(\ker\p)}^{p-1,q+1}=0$ means for any $y\in \ker\p\pb\cap A^{p-1,q+1}(X)$ there is $x\in \ker\p\cap A^{p,q}(X)$ such that $\p y=\pb x$; $\p_{A,BC}^{p,q+1}=0$ means for any $y\in \ker\p\pb\cap A^{p,q+1}(X)$ there is $x\in A^{p-1,q}(X)$ such that $\p y=\p\pb x$. Note that $\p_{A,BC}^{p,q+1}=0$ is exactly the $(p+1,q+1)$-th mild $\p\pb$-lemma introduced by Rao-Wan-Zhao \cite[Def.\,3.1]{RWZ19}. In fact, it was shown in \cite[Coro.\,5.2]{RWZ19} that $\dim H_{BC}^{p,0}(X_t)=\dim H_{A}^{n,n-p}(X_t)$ is independent of $t\in B$ if the $(p,1)$-th mild $\p\pb$-lemma, i.e. $\p_{A,BC}^{p-1,1}=0$ holds. Our Corollary \ref{coro-0} slightly improves this result.

A jumping formula for the Dolbeault cohomology was obtained in \cite[Thm.\,5.10]{Xia19dDol}. By comparing it to the formulas in Theorem \ref{thm-BCandA-jump-formula-0}, it is natural to wonder whether $v^{\bullet,\bullet}_t$ is somehow related to $w^{\bullet,\bullet}_t$. In section \ref{sec-examples}, we present two examples which shows it is possible that $w^{\bullet,\bullet}_t\equiv0$ while at the same time there are some $(p,q)$ with $v^{p,q}_t\neq0$, see Example \ref{example Case III-(2)}.

In \cite{PSU21}, Popovici-Stelzig-Ugarte have introduced the higher-page Bott-Chern and Aeppli cohomologies, it would be interesting to study the corresponding deformation theory and jumping formulas of these higher-page cohomologies.

\vskip 1\baselineskip \textbf{Acknowledgements.} We would like to thank Prof. Daniele Angella, Hisashi Kasuya and Xuanming Ye for useful communications.

\section{Canonical representatives of the deformed Aeppli cohomology}
In this section we study canonical representatives of the deformed Aeppli cohomology. We show for any deformed Aeppli cohomology class $\alpha\in H_{A\phi(t)}^{p,q}(X)$ there is a (possibly non-unique) canonical representative $y\in\alpha$ with $d^*y=0$. If the harmonic part of $\alpha$ is fixed, then such $y$ is uniquely determined and can be expressed as a power series (Proposition \ref{prop-unique small solution}). This fact will be important for our later study about the relations between jumping of $\dim H_{A\phi(t)}^{p,q}(X)$ and the deformation behavior of forms.
\subsection{Hodge theory for Bott-Chern and Aeppli cohomology}
We first recall some basic facts about Hodge theory which will be used later (see \cite{Sch07,Pop15} for more detailed discussions). Let $X$ be a compact complex manifold equipped with a fixed Hermitian metric, then we have the Hodge star operator
\[
*:A^{p,q}(X)\to A^{n-q,n-p}(X).
\]
The Bott-Chern/Aeppli cohomology of $X$ is defined by
\begin{equation}
H^{p,q}_{BC}(X) := \frac{\ker d \cap A^{p,q}(X)}{\im \p\pb \cap A^{p,q}(X)},\quad H^{p,q}_{A}(X) := \frac{\ker \p\pb \cap A^{p,q}(X)}{(\im \p+ \im \pb) \cap A^{p,q}(X)}.
\end{equation}
\subsubsection{Bott-Chern/Aeppli Laplacians and Hodge decompositions}
The Bott-Chern/Aeppli Laplacian operator of $X$ is defined respectively as
\begin{equation}
\begin{split}
\square_{BC}:=&(\p\pb)(\p\pb)^*+(\p\pb)^*(\p\pb)+ (\pb^*\p)(\pb^*\p)^*+(\pb^*\p)^*(\pb^*\p)+ \pb^*\pb + \p^*\p,\\
\square_{A}:=&(\p\pb)^*(\p\pb)+(\p\pb)(\p\pb)^*+ (\pb\p^*)^*(\pb\p^*)+(\pb\p^*)(\pb\p^*)^*+ \pb\pb^* + \p\p^*,
\end{split}
\end{equation}
where $\p^*=-*\pb*$, $\pb^*=-*\p*$ and $(\p\pb)^*=(-1)^{p+q+1}*\p\pb*$ on $A^{p,q}(X)$. Set
\begin{equation}
\begin{split}
\mathcal{H}^{p,q}_{BC}(X):&=\ker\square_{BC}\cap A^{p,q}(X)=\ker d\cap\ker (\p\pb)^*\cap A^{p,q}(X),\\
\mathcal{H}^{p,q}_{A}(X):&=\ker\square_{A}\cap A^{p,q}(X)=\ker\p\pb\cap\ker d^*\cap A^{p,q}(X).
\end{split}
\end{equation}
Then the Hodge star operator induces an isomorphism
\[
*:H^{p,q}_{BC}(X)\cong\mathcal{H}^{p,q}_{BC}(X)\longrightarrow\mathcal{H}^{n-q,n-p}_{A}(X)\cong H^{n-q,n-p}_{A}(X).
\]
Since $\square_{BC}, \square_{A}$ are elliptic, the following orthogonal direct sum decompositions hold:
\begin{equation}\label{eq-BC-Hodgedecomposition}
\begin{split}
A^{\bullet,\bullet}(X)&=\ker\square_{BC}\oplus \im\p\pb\oplus (\im\p^*+\im\pb^*)\\
&=\ker\square_{A}\oplus \im(\p\pb)^*\oplus (\im\p+\im\pb),
\end{split}
\end{equation}
which is equivalent to the existence of the (Bott-Chern/Aeppli) Green operator $G_{BC}, G_{A}$, such that
\begin{equation}\label{eq-BC-Hodgedecomposition-green}
1=\mathcal{H}_{BC}+\square_{BC}G_{BC}=\mathcal{H}_{A}+\square_{A}G_{A},
\end{equation}
where $\mathcal{H}_{BC}, \mathcal{H}_{A}$ are linear projections onto corresponding harmonic spaces.
\subsubsection{Deformed Bott-Chern/Aeppli Laplacians and Hodge decompositions}
Now, let $\pi: (\mathcal{X}, X)\to (B,0)$ be a complex analytic family with Beltrami differentials $\phi(t)\in A^{0,1}(X,T^{1,0})$. We have the deformed version of Bott-Chern/Aeppli Laplacian operators
\begin{equation}
\begin{split}
\square_{BC\phi(t)}:=&[\p\pb_{\phi(t)},(\p\pb_{\phi(t)})^*]+ [\pb_{\phi(t)}^*\p,(\pb_{\phi(t)}^*\p)^*]+ \pb_{\phi(t)}^*\pb_{\phi(t)} + \p^*\p,\\
\square_{A\phi(t)}:=&[\p\pb_{\phi(t)},(\p\pb_{\phi(t)})^*]+ [(\pb_{\phi(t)}\p^*)^*,\pb_{\phi(t)}\p^*]+ \pb_{\phi(t)}\pb_{\phi(t)}^* + \p\p^*.
\end{split}
\end{equation}
Since $\square_{BC\phi(t)}, \square_{A\phi(t)}$ are still elliptic, deformed versions of \eqref{eq-BC-Hodgedecomposition} and \eqref{eq-BC-Hodgedecomposition-green} also holds:
\begin{equation}\label{eq-BC-Hodgedecomposition-phit}
\begin{split}
A^{\bullet,\bullet}(X)&=\ker\square_{BC\phi(t)}\oplus \im\p\pb_{\phi(t)}\oplus (\im\p^*+\im\pb_{\phi(t)}^*)\\
&=\ker\square_{A\phi(t)}\oplus \im(\p\pb_{\phi(t)})^*\oplus (\im\p+\im\pb_{\phi(t)}),
\end{split}
\end{equation}
Set
\begin{equation}
\begin{split}
\mathcal{H}^{p,q}_{BC\phi(t)}(X):&=\ker\square_{BC\phi(t)}\cap A^{p,q}(X)=\ker d_{\phi(t)}\cap\ker (\p\pb_{\phi(t)})^*\cap A^{p,q}(X),\\
\mathcal{H}^{p,q}_{A\phi(t)}(X):&=\ker\square_{A\phi(t)}\cap A^{p,q}(X)=\ker\p\pb_{\phi(t)}\cap\ker d^*_{\phi(t)}\cap A^{p,q}(X).
\end{split}
\end{equation}
\subsubsection{Formal adjoint of $\p\pb_{\phi}$}
\begin{proposition}\label{prop-formal-adjoint}
For any $\phi\in A^{0,1}(X, T^{1,0})$, let $i_{\phi}^*, (\mathcal{L}_{\phi}^{1,0})^*, \pb_{\phi}^*, (\p\pb_{\phi})^*$ be the formal adjoint of $i_{\phi}, \mathcal{L}_{\phi}^{1,0}, \pb_{\phi}=\pb-\mathcal{L}_{\phi}^{1,0}, \p\pb_{\phi}$, respectively. Then the following
\begin{itemize}
\item $i_{\phi}^*=(-1)^{p+q+1}*(i_{\bar{\phi}})*$;
\item $(\mathcal{L}_{\phi}^{1,0})^*=-*(\mathcal{L}_{\bar{\phi}}^{0,1})*$;
\item $\pb_{\phi}^*=-*(\p_{\bar{\phi}})*$;
\item $(\p\pb_{\phi})^*=(-1)^{p+q+1}*(\p_{\bar{\phi}}\pb)*$,
\end{itemize}
holds on $A^{p,q}(X)$ where
\[
\p_{\bar{\phi}}:=\p-\mathcal{L}_{\bar{\phi}}^{0,1}=\p-(i_{\bar{\phi}}\pb-\pb i_{\bar{\phi}}).
\]
\end{proposition}
\begin{proof}Let $\xi\in A^{p+1,q-1}(X)$ and $\eta\in A^{p,q}(X)$, then $0=\xi\wg\overline{*\eta}$. As a result,
\[
0=i_{\phi}(\xi\wg\overline{*\eta})=i_{\phi}\xi\wg\overline{*\eta}+ \xi\wg i_{\phi}\overline{*\eta}\Longrightarrow i_{\phi}\xi\wg\overline{*\eta}=-\xi\wg i_{\phi}\overline{*\eta}=-\xi\wg \overline{**^{-1}i_{\bar{\phi}}*\eta},
\]
which implies
\[
(i_{\phi}\xi,\eta)_{L^2}=\int_X i_{\phi}\xi\wg\overline{*\eta}=-\int_X \xi\wg i_{\phi}\overline{*\eta}=-\int_X \xi\wg \overline{**^{-1}i_{\bar{\phi}}*\eta}=-(\xi,*^{-1}i_{\bar{\phi}}*\eta)_{L^2}.
\]
Hence,
\[
i_{\phi}^*\eta=-*^{-1}i_{\bar{\phi}}*\eta=(-1)^{p+q+1}*i_{\bar{\phi}}*\eta.
\]
For $\mathcal{L}_{\phi}^{1,0}=i_{\phi}\p-\p i_{\phi}$, we have
\begin{align*}
(\mathcal{L}_{\phi}^{1,0})^*=&(i_{\phi}\p-\p i_{\phi})^*\\
=&\p^* i_{\phi}^*-i_{\phi}^*\p^*\\
=&-*\pb*(-1)^{p+q+1}*~i_{\bar{\phi}}~*-(-1)^{p+q}*~i_{\bar{\phi}}~*(-1)*\pb*\\
=&*\pb i_{\bar{\phi}}*-*i_{\bar{\phi}}\pb*=-*\mathcal{L}_{\bar{\phi}}^{0,1}*.
\end{align*}
It follows that
\begin{align*}
\pb_{\phi}^*=&(\pb-\mathcal{L}_{\phi}^{1,0})^*\\
=&\pb^*-(\mathcal{L}_{\phi}^{1,0})^*\\
=&-*\p*~+~*\mathcal{L}_{\bar{\phi}}^{0,1}*=-*(\p_{\bar{\phi}})*,
\end{align*}
and
\[
(\p\pb_{\phi})^*=\pb_{\phi}^*\p^*=*\p_{\bar{\phi}}**\pb*=(-1)^{p+q+1}*(\p_{\bar{\phi}}\pb)*.
\]
\end{proof}

As a direct consequence, we have
\begin{corollary}\label{coro-BCphi-Aphi}
The following homomorphism
\[
\mathcal{H}^{p,q}_{BC\phi(t)}(X)\to\mathcal{H}^{n-p,n-q}_{A\phi(t)}(X):\varphi\longmapsto \overline{*\varphi},
\]
is an isomorphism.
\end{corollary}
\begin{proof}This follows since $*(\ker d_{\phi(t)})=\ker d_{\overline{\phi(t)}}^*$ and $*(\ker (\p\pb_{\phi(t)})^*)=\ker \pb\p_{\overline{\phi(t)}}$.
\end{proof}
Hence, the duality between deformed Bott-Chern and Aeppli cohomology becomes
\begin{equation}\label{eq-BC-A-dual}
H^{p,q}_{BC\phi(t)}(X)\cong H^{n-p,n-q}_{A\phi(t)}(X).
\end{equation}

\subsection{$d^*$-closed representatives of the deformed Aeppli cohomology}
\begin{proposition}\label{prop-unique small solution}The following holds:
\begin{enumerate}
  \item For any $\sigma \in A^{p,q}(X)$, if $\p\pb_{\phi(t)} \sigma=0$ and $\p^*\sigma=\pb^*\sigma=0$, then we must have
\[
\sigma=\mathcal{H}_{A}\sigma + G_{A}(\p\bar{\p})^*\p i_{\phi(t)}\p \sigma,
\]
where $\mathcal{H}_{A}: A^{p,q}(X)\to \mathcal{H}^{p,q}_{A}(X)$ is the projection operator to harmonic space;
  \item For any fixed $\sigma_0\in \mathcal{H}^{p,q}_{A}(X)$ and small $t$, the equation
\begin{equation}\label{Kuranishi eq}
\sigma=\sigma_0 + G_{A}(\p\bar{\p})^*\p i_{\phi(t)}\p \sigma,
\end{equation}
has an unique solution given by $\sigma=\sum_{k} \sigma_k\in A^{p,q}(X)$ and
\[
\sigma_{k}=G_{A}(\p\bar{\p})^*\sum_{i+j=k} \p i_{\phi_j} \p\sigma_i,\quad\forall k > 0;
\]
  \item $\ker\p\pb_{\phi(t)}\cap\im(\p\bar{\p})^*=\ker d_{\phi(t)}\cap\im d^*=\im d\cap\im(\p\bar{\p}_{\phi(t)})^*=0$ for small $t$.
\end{enumerate}
\end{proposition}
\begin{proof}$(1)$ follows from the Hodge decomposition:
\[
\sigma = \mathcal{H}_{A}\sigma + G_{A}\square_{A} \sigma
       = \mathcal{H}_{A}\sigma + G_{A}(\p\bar{\p})^*\p\bar{\p} \sigma
       = \mathcal{H}_{A}\sigma + G_{A}(\p\bar{\p})^*\p i_{\phi(t)}\p \sigma~,
\]
where we have used the fact that $\square_{A} \sigma=(\p\bar{\p})^*\p\bar{\p} \sigma$ if $\p^*\sigma=\pb^*\sigma=0$.\\
$(2)$ Substitute $\sigma=\sigma(t)=\sum_{k}\sigma_k$ in \eqref{Kuranishi eq}, we have
\begin{equation}\label{formal solutions}
\left\{
\begin{array}{ll}
\sigma_1 &= G_{A}(\p\bar{\p})^* \p i_{\phi_1}\p \sigma_0,  \\
\sigma_2 &= G_{A}(\p\bar{\p})^* (\p i_{\phi_2}\p \sigma_0 + \p i_{\phi_1}\p \sigma_1 ),\\
         &\cdots,  \\
\sigma_k &= G_{A}(\p\bar{\p})^*\sum_{i+j=k} \p i_{\phi_j} \p\sigma_i ,~~\forall k > 0.  \\
\end{array} \right.
\end{equation}
For the convergence of $\sigma(t)$, we note that
\begin{equation}
\|\sigma_{j}\|_{k+\alpha}=\| G_{A}(\p\bar{\p})^*\sum_{a+b=j} \p i_{\phi_a}\p \sigma_b \|_{k+\alpha} \leq C \sum_{a+b=j} \|\phi_{a}\|_{k+\alpha} \|\sigma_{b}\|_{k+\alpha},
\end{equation}
for some constant $C$ depends only on $k$ and $\alpha$.
Now it is left to show the uniqueness. Let $\sigma$ and $\sigma'$ be two solutions to $\sigma=\sigma_0+G_{A}(\p\bar{\p})^* \p i_{\phi}\p \sigma$ and set $\tau=\sigma-\sigma'$. Then $\tau=G_{A}(\p\bar{\p})^* \p i_{\phi}\p \tau$, we have
\begin{equation}
\|\tau\|_{k+\alpha}\leq c \|\phi (t)\|_{k+\alpha}  \|\tau\|_{k+\alpha},
\end{equation}
for some constant $c>0$. When $|t|$ is sufficiently small, $\|\phi (t)\|_{k+\alpha}$ is also small. Hence we must have $\tau=0$.
For smoothness of the solution, note that we have
\[
\square_{A}\sigma=(\p\bar{\p})^*\p\bar{\p} \sigma=(\p\bar{\p})^*\p i_{\phi(t)}\p \sigma,
\]
which implies
\begin{equation}\label{regularity eq}
\square_{A}\sigma - (\p\bar{\p})^*\p i_{\phi(t)}\p\sigma= 0,
\end{equation}
which is a standard elliptic equation for small $t$.\\
$(3)$ Assume $\sigma\in\ker\p\pb_{\phi(t)}\cap\im(\p\bar{\p})^*$. Then we have $\mathcal{H}_{A}\sigma=0$ and $\sigma=\mathcal{H}_{A}\sigma + G_{A}(\p\bar{\p})^*\p i_{\phi(t)}\p \sigma$ by $(1)$. It follows from $(2)$ that $\sigma=0$.

Similarly, assume $\sigma\in\ker d_{\phi(t)}\cap\im d^*$. Then
\[
\sigma=\sigma_0 - G_{BC}(\pb^*\p\p^*+\pb^*)\p i_{\phi(t)} \sigma
\]
with $\sigma_0:=\mathcal{H}_{BC}\sigma$. Applying the Bott-Chern harmonic projection $\mathcal{H}_{BC}$ to the both sides of this equation, we immediately get $\sigma_0=0$. Hence we must have $\sigma=0$.

Assume $\tilde{\sigma}\in\im d\cap\im(\p\bar{\p}_{\phi(t)})^*\cap A^{p,q}(X)$. Then it follows from Proposition \ref{prop-formal-adjoint} that
\[
*\tilde{\sigma}\in\im d^*\cap\im(\p_{\overline{\phi(t)}}\pb)\cap A^{n-q,n-p}(X)
\]
which means
\[
\overline{*\tilde{\sigma}}\in\im d^*\cap\im\p\pb_{\phi(t)}\cap A^{n-p,n-q}(X).
\]
So in order to prove $\im d\cap\im(\p\bar{\p}_{\phi(t)})^*=0$ we only need to show $\im d^*\cap\im\p\pb_{\phi(t)}=0$. This is clear since
\[
\im d^*\cap\im\p\pb_{\phi(t)}\subseteq \im d^*\cap\ker d_{\phi(t)}=0.
\]
\end{proof}

\begin{proposition}\label{prop-Aclosed=BCharmonic part vanish}
Let $\sigma=\sigma(t)$ be a solution of the equation \eqref{Kuranishi eq} in Proposition \ref{prop-unique small solution}.
Then for any small $t\in B$, we have
\begin{equation}\label{eq-equivalence-BCharmonic part vanish}
\p\pb_{\phi(t)}\sigma(t) = 0  \Leftrightarrow \mathcal{H}_{BC}\Big(\p i_{\phi(t)}\p\sigma(t)\Big)=0.
\end{equation}
Moreover, we have $\p i_{\phi(t)}\p\sigma(t)\in\ker d$ for any small $t\in B$.
\end{proposition}
\begin{proof}If $\p\pb_{\phi(t)}\sigma(t) = 0$, then $\p\pb\sigma(t)=\p i_{\phi(t)}\p\sigma(t)$ which implies $\mathcal{H}_{BC}\Big(\p i_{\phi(t)}\p\sigma(t)\Big)=0$.

Conversely, suppose $\mathcal{H}_{BC}\Big(\p i_{\phi(t)}\p\sigma(t)\Big)=0$. Then
\begin{align*}
\p i_{\phi(t)}\p\sigma(t)
=&G_{BC}\square_{BC}\p i_{\phi(t)}\p\sigma(t)\\
=&G_{BC}\Big( (\p\pb)(\p\pb)^*+ (\pb^*\p)(\pb^*\p)^*+\pb^*\pb\Big) \p i_{\phi(t)}\p\sigma(t).
\end{align*}
Set
\[
\psi(t):=\p\pb_{\phi(t)}\sigma(t)=\p\pb\sigma(t)-\p i_{\phi(t)}\p\sigma(t).
\]
Recall that since the Maurer-Cartan equation $\pb\phi(t)=\frac{1}{2}[\phi(t),\phi(t)]$ holds, we have
\begin{equation}\label{eq-liepbphi}
\mathcal{L}^{1,0}_{\pb\phi(t)}=\frac{1}{2}\mathcal{L}^{1,0}_{[\phi(t),\phi(t)]}=\frac{1}{2}[\mathcal{L}^{1,0}_{\phi(t)},\mathcal{L}^{1,0}_{\phi(t)}]=\mathcal{L}^{1,0}_{\phi(t)}\mathcal{L}^{1,0}_{\phi(t)},
\end{equation}
By substituting \eqref{Kuranishi eq} and using \eqref{eq-liepbphi}, we get
\begin{align*}
\psi(t)=&\p\pb\Big(G_{A}(\p\bar{\p})^*\p i_{\phi(t)}\p \sigma(t)\Big)-\p i_{\phi(t)}\p\sigma(t)\\
=&G_{BC}\p\pb(\p\bar{\p})^*\p i_{\phi(t)}\p \sigma(t)-\p i_{\phi(t)}\p\sigma(t)\\
=&-G_{BC}\Big( \pb^*\p\p^*+\pb^* \Big)\pb\p i_{\phi(t)}\p\sigma(t)\\
=&G_{BC}\Big( \pb^*\p\p^*+\pb^* \Big)\p\pb\mathcal{L}^{1,0}_{\phi(t)}\sigma(t)\\
=&G_{BC}\Big( \pb^*\p\p^*+\pb^* \Big)\p(-\mathcal{L}^{1,0}_{\phi(t)}\pb+\mathcal{L}^{1,0}_{\pb\phi(t)} )\sigma(t)\\
=&G_{BC}\Big( \pb^*\p\p^*+\pb^* \Big)\p(-i_{\phi(t)}\p\pb+\mathcal{L}^{1,0}_{\phi(t)}\mathcal{L}^{1,0}_{\phi(t)} )\sigma(t)\\
=&G_{BC}\Big( \pb^*\p\p^*+\pb^* \Big)\p i_{\phi(t)}(-\p\pb+\p i_{\phi(t)}\p )\sigma(t)\\
=&-G_{BC}\Big( \pb^*\p\p^*+\pb^* \Big)\p i_{\phi(t)}\psi(t)\\
\end{align*}
which implies
\[
\|\psi(t)\|_{k+\alpha}=\| G_{BC}\Big( \pb^*\p\p^*+\pb^* \Big)\p i_{\phi(t)}\psi(t) \|_{k+\alpha} \leq C \|\phi(t)\|_{k+\alpha}\|\psi(t)\|_{k+\alpha},
\]
where $C$ depends only on $k$ and $\alpha$. Now choose small $|t|$ so that $C \|\phi(t)\|_{k+\alpha}<1$, then we get $\|\psi(t)\|_{k+\alpha}< \|\psi(t)\|_{k+\alpha}$ which is a contradiction if $\psi(t)\neq0$. So we must have $\psi(t)=0$ whenever $|t|$ is small enough. This proves \eqref{eq-equivalence-BCharmonic part vanish}.

Finally, since by \eqref{Kuranishi eq} we have $\p\pb\sigma(t)=\p\pb G_{A}(\p\bar{\p})^*\p i_{\phi(t)}\p \sigma$, it follows that
\begin{align*}
&\p i_{\phi(t)}\p\sigma(t)\\
=&\mathcal{H}_{BC}\Big(\p i_{\phi(t)}\p\sigma(t)\Big)+G_{BC}\Big( (\p\pb)(\p\pb)^*+ (\pb^*\p)(\pb^*\p)^*+\pb^*\pb\Big) \p i_{\phi(t)}\p\sigma(t)\\
=&\mathcal{H}_{BC}\Big(\p i_{\phi(t)}\p\sigma(t)\Big)+\Big((\p\pb)G_{A}(\p\pb)^*+G_{BC}(\pb^*\p)(\pb^*\p)^*+G_{BC}\pb^*\pb\Big) \p i_{\phi(t)}\p\sigma(t)\\
=&\mathcal{H}_{BC}\Big(\p i_{\phi(t)}\p\sigma(t)\Big)+\p\pb\sigma(t)+G_{BC}\Big((\pb^*\p)(\pb^*\p)^*+\pb^*\pb\Big) \p i_{\phi(t)}\p\sigma(t),
\end{align*}
which combined with \eqref{eq-equivalence-BCharmonic part vanish} will imply
\[
G_{BC}\Big((\pb^*\p)(\pb^*\p)^*+\pb^*\pb\Big) \p i_{\phi(t)}\p\sigma(t)=0.
\]
As a result, we get
\[
\p i_{\phi(t)}\p\sigma(t)\\
=\mathcal{H}_{BC}\Big(\p i_{\phi(t)}\p\sigma(t)\Big)+G_{BC}(\p\pb)(\p\pb)^*\p i_{\phi(t)}\p\sigma(t)\in\ker d.
\]
\end{proof}
\begin{remark}In the Bott-Chern case, similar proofs will show that if $\sigma=\sigma(t)$ is a solution to $\sigma=\sigma_0 - G_{BC}(\pb^*\p\p^*+\pb^*) \p i_{\phi(t)} \sigma$ for some given $\sigma_0\in \mathcal{H}^{p,q}_{BC}(X)$, then

\end{remark}

For a vector space $A$ and its subspaces $V,W,U$, if we have a direct sum decomposition
\[
A=V\oplus W,
\]
then the following homomorphism
\begin{equation}\label{eq-U=Vcap(U+W)}
\frac{U}{U\cap W}\longrightarrow V\cap(U+W):[u]\longmapsto v,
\end{equation}
where $u=v+w\in V\oplus W$, is an isomorphism. In case that $U\cap W=\{0\}$, then $U\cong V\cap(U+W)$. This simple observation will be used frequently.
\begin{proposition}\label{prop-d*-closed-rep of deformed-A-coho}
The natural map
\begin{equation}\label{eq-natural map}
\frac{\ker d^*\cap\ker \p\pb_{\phi(t)} \cap A^{p,q}(X)}{\ker d^*\cap (\im \p+ \im \pb_{\phi(t)}) \cap A^{p,q}(X)}\longrightarrow H_{A\phi(t)}^{p,q}(X)
\end{equation}
is an isomorphism.
\end{proposition}
\begin{proof}By \eqref{eq-BC-Hodgedecomposition} and \eqref{eq-BC-Hodgedecomposition-phit}, we have the following orthogonal direct sum decomposition\footnote{This is just the Hodge decomposition for the Laplacian\[\square_{A}^{\phi(t)}:=(\p\pb_{\phi(t)})^*(\p\pb_{\phi(t)})+(\p\pb)(\p\pb)^*+ (\pb\p^*)^*(\pb\p^*)+(\pb\p^*)(\pb\p^*)^*+ \pb\pb^* + \p\p^*.\] }
\[
A^{p,q}(X)=\left(\ker d^*\cap\ker\p\pb_{\phi(t)}\right)\oplus
(\im\p+\im\pb+\im(\p\pb_{\phi(t)})^*)
\]
which implies
\[
\ker\p\pb_{\phi(t)}=\left(\ker d^*\cap\ker\p\pb_{\phi(t)}\right)\oplus
\left(\ker\p\pb_{\phi(t)}\cap(\im\p+\im\pb+\im(\p\pb_{\phi(t)})^*)\right),
\]
and
\begin{align*}
\im\p+\im\pb_{\phi(t)}=&\left(\ker d^*\cap(\im\p+\im\pb_{\phi(t)})\right)\oplus\\
&\left((\im\p+\im\pb_{\phi(t)})\cap(\im\p+\im\pb+\ker\p^*\cap\ker\pb_{\phi(t)}^*)\right).
\end{align*}
Moreover, according to \eqref{eq-U=Vcap(U+W)}, we have
\[
\ker\p\pb_{\phi(t)}\cap\left(\im\p+\im\pb+\im(\p\pb_{\phi(t)})^*\right)\cong\frac{\im\p+\im\pb}{(\im\p+\im\pb)\cap\im(\p\pb_{\phi(t)})^*}
\]
and
\[
(\im\p+\im\pb_{\phi(t)})\cap(\im\p+\im\pb+\ker\p^*\cap\ker\pb_{\phi(t)}^*)
\cong\frac{\im\p+\im\pb}{(\im\p+\im\pb)\cap\ker\p^*\cap\ker\pb_{\phi(t)}^*}.
\]
Hence,
\[
H_{A\phi(t)}^{p,q}(X)\cong
\frac{\ker d^*\cap\ker\p\pb_{\phi(t)} }{\ker d^*\cap(\im\p+\im\pb_{\phi(t)}) } \oplus
\frac{(\im\p+\im\pb)\cap\ker\p^*\cap\ker\pb_{\phi(t)}^* }{(\im\p+\im\pb)\cap\im(\p\pb_{\phi(t)})^* }
\]
The conclusion now follows because $(\im\p+\im\pb)\cap\ker\p^*\cap\ker\pb_{\phi(t)}^*=0$ by $(3)$ of Proposition \ref{prop-unique small solution}. Indeed, note that
\[
(\im\p+\im\pb)\cap\ker\p^*\cap\ker\pb_{\phi(t)}^*=\overline{*(\ker d_{\phi(t)}\cap\im d^*)}.
\]
\end{proof}

\section{Canonical Aeppli deformations of $(p,q)$-forms}
In this section, we introduce the notion of canonical Aeppli deformations of $(p,q)$-forms and study some of its basic properties.
\subsection{The notions of Dolbeault/Bott-Chern/Aeppli deformations of $(p,q)$-forms}
Let $\pi: (\mathcal{X}, X)\to (B,0)$ be a complex analytic family such that for each $t\in B$ the complex structure on $X_t$ is represented by Beltrami differential $\phi(t)$. In order to motivate the definition of Aeppli deformations of $(p,q)$-forms, we will also recall the definitions of Dolbeault/Bott-Chern deformations of $(p,q)$-forms studied in \cite{Xia19dDol,Xia19dBC}.
\begin{definition}\label{def-deformation-Dol/BC/A}
Given a $(p,q)$-form
\[
y\in \ker\pb/\ker d/\ker\p\pb
\]
and $T\subseteq B$, which is an analytic subset of $B$ containing $0$,
\[
\text{a~\emph{Dolbeault/Bott-Chern/Aeppli deformation}~of~$y$~w.r.t.}~\pi: (\mathcal{X}, X)\to (B,0)~\text{on}~T
\]
is a family of $(p,q)$-forms $\sigma (t)$ such that
\begin{itemize}
  \item[1.] $\sigma (t)$ is holomorphic in $t\in T$ and $\sigma (0) = y$;
  \item[2.] $\forall t\in T$, the following holds respectively:
  \begin{itemize} \item $\pb_{\phi(t)}\sigma (t) = 0$ (Dolbeault);
  \item $d_{\phi(t)}\sigma (t) = 0$ (Bott-Chern);
  \item $\p\pb_{\phi(t)}\sigma (t) = 0$ (Aeppli);
  \end{itemize}
\end{itemize}
A \emph{Dolbeault/Bott-Chern/Aeppli deformation} of a class
\[
\alpha\in H_{\pb}^{p,q}(X)/H_{BC}^{p,q}(X)/H_{A}^{p,q}(X)
\]
(w.r.t. $\pi$) on $T$ is a triple $(y,\sigma (t),T)$ which consisting of a representative $y\in\alpha$ and a Dolbeault/Bott-Chern/Aeppli deformation $\sigma (t)$ of $y$ (w.r.t. $\pi$) on $T$. Two deformations $(y,\sigma (t),T)$ and $(y',\sigma' (t),T)$ of $[y]$ of $\alpha$ on $T$ are \emph{equivalent} if
\[
[\sigma (t) - \sigma' (t)] = 0 \in H_{\pb_{\phi(t)}}^{p,q}(X)/H_{BC\phi(t)}^{p,q}(X)/H_{A\phi(t)}^{p,q}(X),~\forall t\in T.
\]
\end{definition}
\subsubsection{Canonical deformations}
\begin{definition}\label{def-canonical-deformation}
Assume $X$ is equipped with a fixed Hermitian metric $h$. Correspondingly, we have the following notions for canonical deformations:
\begin{description}
\item[Dolbeault]
A Dolbeault deformation $\sigma (t)$ of $y\in\ker\pb$ on $T$ is called \emph{canonical} if $\sigma (t)$ is of the form
\[
\sigma(t)=\sigma^h(t)+\pb_{\phi(t)}u,
\]
such that
\begin{itemize}
\item $\sigma^h(t) = \mathcal{H}_{\pb}y + \pb^*G \mathcal{L}_{\phi(t)}^{1,0}\sigma^h(t),~\forall t\in T$;
\item $\pb u=y-\mathcal{H}_{\pb}y$.
\end{itemize}
\item[Bott-Chern] A Bott-Chern deformation $\sigma (t)$ of $y\in\ker d$ on $T$ is called \emph{canonical} if $\sigma (t)$ is of the form
\[
\sigma(t)=\sigma^h(t)+\p\pb_{\phi(t)}u,
\]
such that
\begin{itemize}
\item $\sigma^h(t) = \mathcal{H}_{BC}y - G_{BC}(\pb^*\p\p^*+\pb^*) \p i_{\phi(t)} \sigma^h(t),~\forall t\in T$;
\item $\p\pb u=y-\mathcal{H}_{BC}y$.
\end{itemize}
\item[Aeppli] An Aeppli deformation $\sigma (t)$ of $y\in\ker\p\pb$ on $T$ is called \emph{canonical} if $\sigma (t)$ is of the form
\[
\sigma(t)=\sigma^h(t)+\p v+\pb_{\phi(t)}u,
\]
such that
\begin{itemize}
\item $\sigma^h(t) = \mathcal{H}_{A}y + G_{A}(\p\bar{\p})^*\p i_{\phi(t)}\p \sigma^h(t),~\forall t\in T$;
\item $\p v+\pb u=y-\mathcal{H}_{A}y$.
\end{itemize}
\end{description}
\end{definition}
Under this definition, Dolbeault/Bott-Chern/Aeppli-exact forms can always deform canonically into a family of deformed Dolbeault/Bott-Chern/Aeppli-exact forms. For example, in the Aeppli case, if $y=\p v+\pb u$, then $\sigma(t)=\p v+\pb_{\phi(t)}u$ is clearly a canonical (Aeppli) deformation of $y$ on $B$. In fact, we have
\begin{proposition}\label{prop-unique-can}
Let $\alpha\in H_{\pb}^{p,q}(X)/H_{BC}^{p,q}(X)/H_{A}^{p,q}(X)$. Suppose a Hermitian metric on $X$ is fixed. Then any two canonical Dolbeault/Bott-Chern/Aeppli deformation of $\alpha$ (on some common domain $T$) are equivalent.
\end{proposition}
\begin{proof}Let us take the case of canonical Aeppli deformation as an example, the other two cases are similar. Indeed, let $y, y'\in\alpha$ and $\alpha\in H_{A}^{p,q}(X)$.
If $\sigma(t)=\sigma^h(t)+\p v+\pb_{\phi(t)}u$ and $\sigma'(t)=\sigma^{'h}(t)+\p v'+\pb_{\phi(t)}u'$ are canonical Aeppli deformation of $y, y'$ (on $T$), respectively. Then by Proposition \ref{prop-unique small solution}, we have $\sigma^h(t)=\sigma^{'h}(t)$ since $\mathcal{H}_{A}y=\mathcal{H}_{A}y'$. Now,
\[
\sigma(t)-\sigma'(t)=\p(v-v')+\pb_{\phi(t)}(u-u')\in\im\p+\im\pb_{\phi(t)},
\]
from which we conclude that $\sigma(t)$ is equivalent to $\sigma'(t)$.
\end{proof}
\begin{remark}In our previous work \cite{Xia19dDol,Xia19dBC}, canonical deformations have slightly different meanings. They coincide with Definition \ref{def-canonical-deformation} for (Dolbeault/Bott-Chern) harmonic forms. In view of the principle ``exact forms/zero classes can always deforms into a family of deformed exact forms/zero classes", the definitions given as above seems better suited for our purpose. In particular, it follows immediately from Definition \ref{def-canonical-deformation} that canonical Dolbeault/Bott-Chern/Aeppli deformations of exact forms always exist and are equivalent to the trivial deformation. This is exactly the content of \cite[Lem.\,5.6]{Xia19dBC} (in the Dolbeault case). On the other hand, we point out that under Definition \ref{def-canonical-deformation} canonical deformation of a given closed form is not unique in general. However, by Proposition \ref{prop-unique-can} the equivalence classes of canonical deformations of a given closed form/class are unique as long as we fix the Hermitian metric on $X$.
\end{remark}
\subsubsection{Unobstructedness and canonically unobstructedness}
\begin{definition}\label{def-unobstructedness}
Suppose there is a complex analytic family $\pi: (\mathcal{X}, X)\to (B,0)$.
\begin{enumerate}
  \item A $(p,q)$-form $y\in\ker\pb/\ker d/\ker\p\pb$ is said to be \emph{(canonically) unobstructed w.r.t. $\pi$} if a (canonical) deformation of $y$ (w.r.t. $\pi$) exists on $B$;
  \item A class $\alpha\in H_{\pb}^{p,q}(X)/H_{BC}^{p,q}(X)/H_{A}^{p,q}(X)$ is \emph{(canonically) unobstructed w.r.t. $\pi$} if there is a $y\in\alpha$ such that $y$ is (canonically) unobstructed w.r.t. $\pi$;
  \item If every class in $H_{\pb}^{p,q}(X)/H_{BC}^{p,q}(X)/H_{A}^{p,q}(X)$ have canonically unobstructed deformation w.r.t. $\pi$, then we say \emph{classes in $H_{\pb}^{p,q}(X)/H_{BC}^{p,q}(X)/H_{A}^{p,q}(X)$ are canonically unobstructed w.r.t. $\pi$} or \emph{the Dolbeault/Bott-Chern/Aeppli deformations of $(p,q)$-forms on $X$ is canonically unobstructed w.r.t. $\pi$}. If this holds for any other complex analytic family of $X$, we will accordingly drop the term ``w.r.t. $\pi$''.
\end{enumerate}
\end{definition}
For any $y\in\mathcal{H}_{A}^{p,q}(X)$, if $p=n$ or $q=n$ then it is always canonically unobstructed by definition: its canonical deformation is given by $\sigma(t)\equiv y$.

\subsection{The dimension of $\ker(\p\pb)^*\cap\im\p\pb_{\phi(t)}$}
Recall the following:
\begin{proposition}\cite[Prop.\,4.3]{Xia19dBC}\label{prop-kerppb*-imagep+ppbphi} Let $\pi: (\mathcal{X}, X)\to (B,0)$ be a complex analytic family such that for each $t\in B$ the complex structure on $X_t$ is represented by Beltrami differential $\phi(t)$.
\begin{enumerate}
\item For any fixed $t\in B$, the following homomorphism
\begin{align*}
g_t:\ker\p\pb\cap A^{p,q}(X)&\longrightarrow \ker(\p\pb)^*\cap\im\p\pb_{\phi(t)}\cap A^{p+1,q+1}(X)\\
x_0&\longmapsto \p\pb_{\phi(t)}x(t),
\end{align*}
is surjective with
\[
\ker g_t=\ker\p\pb\cap\left(\ker\p\pb_{\phi(t)}+\im(\p\pb)^*\right)\cap A^{p,q}(X),
\]
where $x(t)$ is the unique solution of $x(t)=x_0 + (\p\pb)^*G_{BC} \p i_{\phi(t)} \p x(t)$.\\
\item Let
\[
\hat{g}_t:\mathcal{H}^{p,q}_{BC}(X)\longrightarrow \ker(\p\pb)^*\cap\im\p\pb_{\phi(t)}\cap A^{p+1,q+1}(X)
\]
be the restriction of $g_t$ on $\mathcal{H}^{p,q}_{BC}(X)$. Then $\hat{g}_t$ is surjective with
\[
\ker \hat{g}_t=\mathcal{H}^{p,q}_{BC}(X)\cap\left(\ker\p\pb_{\phi(t)}+\im(\p\pb)^*\right).
\]
Moreover, we have
\begin{align*}
\dim\mathcal{H}^{p,q}_{BC}(X)=&\dim\ker\p\pb_{\phi(t)}\cap(\mathcal{H}^{p,q}_{BC}(X)+\im(\p\pb)^*)\cap A^{p,q}(X)\\+
&\dim \ker(\p\pb)^*\cap\im\p\pb_{\phi(t)}\cap A^{p+1,q+1}(X).
\end{align*}
\end{enumerate}
\end{proposition}

\begin{proposition}\label{prop-kerd*-imagep+pbarphi}
In the context of Proposition \ref{prop-kerppb*-imagep+ppbphi}, let
\[
\check{g}_t:\mathcal{H}^{p,q}_{A}(X)\longrightarrow \ker(\p\pb)^*\cap\im\p\pb_{\phi(t)}\cap A^{p+1,q+1}(X)
\]
be the restriction of $g_t$ on $\mathcal{H}^{p,q}_{A}(X)$. Then $\check{g}_t$ is surjective with
\[
\ker\check{g}_t=\mathcal{H}^{p,q}_{A}(X)\cap\left(\ker\p\pb_{\phi(t)}+\im(\p\pb)^*\right).
\]
Moreover, we have
\begin{align*}
\dim\mathcal{H}^{p,q}_{A}(X)=&\dim\ker\p^*\cap\ker\pb^*\cap\ker\p\pb_{\phi(t)}\cap A^{p,q}(X)\\+
&\dim \ker(\p\pb)^*\cap\im\p\pb_{\phi(t)}\cap A^{p+1,q+1}(X).
\end{align*}
\end{proposition}
\begin{proof}It follows immediately from $\ker g_t=\ker\p\pb\cap\left(\ker\p\pb_{\phi(t)}+\im(\p\pb)^*\right)$ that
\[
\ker\check{g}_t=\mathcal{H}^{p,q}_{A}(X)\cap\left(\ker\p\pb_{\phi(t)}+\im(\p\pb)^*\right).
\]
To show $\check{g}_t$ is surjective it is enough to show
\[
\frac{\mathcal{H}^{p,q}_{A}(X)}{\mathcal{H}^{p,q}_{A}(X)\cap\left(\ker\p\pb_{\phi(t)}+\im(\p\pb)^*\right)}\cong
\frac{\ker\p\pb\cap A^{p,q}(X)}{\ker\p\pb\cap\left(\ker\p\pb_{\phi(t)}+\im(\p\pb)^*\right)\cap A^{p,q}(X)}.
\]
Indeed, it follows from (see $(3)$ of Proposition \ref{prop-unique small solution})
\[
\im d\cap\im(\p\bar{\p}_{\phi(t)})^*=\ker\p\pb_{\phi(t)}\cap\im(\p\pb)^*=0
\]
that
\begin{align*}
\ker\p\pb=& \mathcal{H}^{p,q}_{A}(X)\oplus(\im\p+ \im\pb)\\
\cong& \mathcal{H}^{p,q}_{A}(X)\oplus[\ker\p\pb_{\phi(t)}\cap(\im\p+\im\pb+\im(\p\pb_{\phi(t)})^*)]
\end{align*}
and
\begin{align*}
&\ker\p\pb\cap\left(\ker\p\pb_{\phi(t)}+\im(\p\pb)^*\right)\cong\ker\p\pb_{\phi(t)}\\
=&\left\{\ker\p\pb_{\phi(t)}\cap (\im\p+\im\pb+\im(\p\pb_{\phi(t)})^*)\right\}\oplus\\
&\left\{\ker\p\pb_{\phi(t)}\cap[\im(\p\pb_{\phi(t)})^*+\ker\p^*\cap\ker\pb^*\cap\ker\p\pb_{\phi(t)}]\right\}\\
\cong& \left\{\ker\p\pb_{\phi(t)}\cap (\im\p+\im\pb+\im(\p\pb_{\phi(t)})^*)\right\}\oplus
\left(\ker\p^*\cap\ker\pb^*\cap\ker\p\pb_{\phi(t)}\right).
\end{align*}
On the other hand, we note that
\begin{align*}
\mathcal{H}^{p,q}_{A}(X)\cap\left(\ker\p\pb_{\phi(t)}+\im(\p\pb)^*\right)=&\mathcal{H}^{p,q}_{A}(X)\cap\left(\ker d^*\cap\ker\p\pb_{\phi(t)}+\im(\p\pb)^*\right)\\
(\text{since}~\ker d^*=\mathcal{H}^{p,q}_{A}(X)\oplus\im(\p\pb)^*)\cong&\ker d^*\cap\ker\p\pb_{\phi(t)}\\
=&\ker\p^*\cap\ker\pb^*\cap\ker\p\pb_{\phi(t)}.
\end{align*}
\end{proof}
\section{The jumping formulas for the deformed Bott-Chern and Aeppli cohomology}
Let $\pi: (\mathcal{X}, X)\to (B,0)$ be a complex analytic family such that for each $t\in B$ the complex structure on $X_t$ is represented by Beltrami differential $\phi(t)$. In this section, we prove jumping formulas for the deformed Bott-Chern and Aeppli cohomology. These formulas precisely relate the canonical Bott-Chern/Aeppli deformations of forms to the jumping behavior of $\dim H_{BC\phi(t)}^{p,q}(X)$, $\dim H_{A\phi(t)}^{p,q}(X)$ with respect to parameter $t$.
\subsection{Jumping of cohomology dimensions and domain of canonical deformations}
For any $t\in B$ and a vector subspace $V\subseteq \mathcal{H}_{BC}^{p,q}(X)$, the following was introduced in \cite[Def.\,4.5]{Xia19dBC}:
\begin{equation}\label{eq-V BC,t}
\begin{split}
V_{BC,t}^{p,q}:=
&\Big\{ \sigma_0\in V \mid  \sigma(t)\in\ker d_{\phi(t)},~\text{where}~\sigma (t)=\sum_{k\geq0} \sigma_k~\text{with}\\
&\sigma_{k}=-G_{BC}(\pb^*\p\p^*+\pb^*)\sum_{i+j=k} \p i_{\phi_j} \sigma_i,~\forall k> 0 \Big\}.
\end{split}
\end{equation}
When $V=\mathcal{H}_{BC}^{p,q}(X)$, we have (c.f. \cite[Prop.\,4.7]{Xia19dBC})
\begin{equation}\label{eq-kertoV}
V_{BC,t}^{p,q}\cong \ker(\p\pb)^*\cap\ker d_{\phi(t)} \cap A^{p,q}(X),
\end{equation}
which implies
\begin{equation}\label{eq-vtoV}
\begin{split}
v^{p,q}_t:=&\dim H_{BC}^{p,q}(X)-\dim \ker d_{\phi(t)}\cap\ker(\p\bar{\p})^*\cap A^{p,q}(X)\\
=&\dim H_{BC}^{p,q}(X)-\dim V_{BC,t}^{p,q}.
\end{split}
\end{equation}
These constructions for Bott-Chern deformations have Aeppli analogues as follows.
\begin{definition}\label{def-W_tp,q-f_t}
For any $t\in B$ and a vector subspace $W\subseteq \mathcal{H}_{A}^{p,q}(X)$, we set
\begin{align*}
W_{t}^{p,q}:=
&\Big\{ \sigma_0\in W \mid  \sigma(t)\in\ker \p\pb_{\phi(t)},\\
&\text{where}~\sigma (t)=\sum_{k} \sigma_k~\text{with}~\sigma_k=G_{A}(\p\bar{\p})^*\sum_{i+j=k} \p i_{\phi_j} \p\sigma_i,~\forall k>0 \Big\},
\end{align*}
and
\begin{align*}
f_t^{p,q}: &W_{t}^{p,q} \longrightarrow \frac{\ker d^*\cap\ker \p\pb_{\phi(t)} \cap A^{p,q}(X)}{\ker d^*\cap (\im \p+ \im \pb_{\phi(t)}) \cap A^{p,q}(X)}\cong H_{A\phi(t)}^{p,q}(X),\\
&\sigma_0\longmapsto \sigma(t)=\sum_{k} \sigma_k,~\text{where}~\sigma_{k}=G_{A}(\p\bar{\p})^*\sum_{i+j=k} \p i_{\phi_j} \p\sigma_i,~\forall k>0.
\end{align*}
\end{definition}
\begin{remark}In view of Proposition \ref{prop-unique-can}, the requirement of $W$ being a subspace of the harmonic space $\mathcal{H}_{A}^{p,q}(X)$ may be removed. More precisely, for any subspace
\[
W_0=\C\{[\sigma_0^1],\cdots,[\sigma_0^l]\}\subset H_{A}^{p,q}(X),
\]
of dimension $l$. Set $\tilde{W}:=\C\{\sigma_0^1,\cdots,\sigma_0^l\}$ and
\[
\tilde{W}_{t}^{p,q}:=
\Big\{ \sigma_0\in \tilde{W} \mid  \sigma(t)\in\ker \p\pb_{\phi(t)},
\text{where}~\sigma (t)~\text{is~a~canonical~deformation~of}~\sigma_0\},
\]
Then it follows from Proposition \ref{prop-unique-can} that we have an isomorphism
\begin{align*}
\mathcal{H}_{A}: &\tilde{W}_{t}^{p,q} \longrightarrow W_{t}^{p,q},\\
&\sigma_0\longmapsto \mathcal{H}_{A}\sigma_0,
\end{align*}
where $W_{t}^{p,q}$ is as in Definition \ref{def-W_tp,q-f_t} for $W=\C\{\mathcal{H}_{A}\sigma_0^1,\cdots,\mathcal{H}_{A}\sigma_0^l\}$.
\end{remark}
By definition, the space $W_{t}^{p,q}$.According to Proposition \ref{prop-unique small solution}, the following linear mapping
\begin{align*}
\tilde{f}_t^{p,q}: &W_{t}^{p,q} \longrightarrow \ker d^*\cap\ker \p\pb_{\phi(t)} \cap A^{p,q}(X),\\
&\sigma_0\longmapsto \sigma(t)=\sum_{k} \sigma_k,~\text{where}~\sigma_{k}=G_{A}(\p\bar{\p})^*\sum_{i+j=k} \p i_{\phi_j} \p\sigma_i,~\forall k\neq 0,
\end{align*}
is an isomorphism for any small $t\in B$. In particular, when $W=\mathcal{H}_{A}^{p,q}(X)$,
\begin{equation}\label{eq-wt-Wt}
\begin{split}
w_{t}^{p,q}:=&\dim H_{A}^{p,q}(X)-\dim\ker\p\pb_{\phi(t)}\cap\ker d^*\cap A^{p,q}(X)\\
=&\dim H_{A}^{p,q}(X) - \dim W_{t}^{p,q},
\end{split}
\end{equation}
and $f_t^{p,q}$ is surjective with $\ker f_t^{p,q}\cong \ker d^*\cap (\im \p+ \im \pb_{\phi(t)})\cap A^{p,q}(X)$. On the other hand, we have
\subsection{The jumping formulas}
\begin{theorem}\label{thm-BCandA-jump-formula}
Let $\pi: (\mathcal{X}, X)\to (B,0)$ be a complex analytic family such that for each $t\in B$ the complex structure on $X_t$ is represented by Beltrami differential $\phi(t)$. Suppose $X$ is equipped with a fixed Hermitian metric.
Then for each $(p,q)\in \mathbb{N}\times \mathbb{N}$, we have the jumping formula for the deformed Bott-Chern cohomology:
\begin{equation}\label{eq-BCjump}
\dim H_{BC}^{p,q}(X)=\dim H_{BC\phi(t)}^{p,q}(X)+v^{p,q}_t+w^{p-1,q-1}_t,
\end{equation}
and the jumping formula for the deformed Aeppli cohomology:
\begin{equation}\label{eq-Ajump}
\dim H_{A}^{p,q}(X)=\dim H_{A\phi(t)}^{p,q}(X)+v^{n-p,n-q}_t+w^{n-p-1,n-q-1}_t.
\end{equation}
\end{theorem}
\begin{proof}First, we know that
\[
H_{BC\phi(t)}^{p,q}(X)\cong\frac{\ker(\p\pb)^*\cap\ker d_{\phi(t)} \cap A^{p,q}(X)}{\ker(\p\pb)^*\cap\im\p\pb_{\phi(t)} \cap A^{p,q}(X)}.
\]
It follows from this, \eqref{eq-vtoV} and Proposition \ref{prop-kerd*-imagep+pbarphi} that
\begin{align*}
\dim H_{BC\phi(t)}^{p,q}(X)=&\dim H_{BC}^{p,q}(X)-v^{p,q}_t-\dim\ker(\p\pb)^*\cap\im\p\pb_{\phi(t)} \cap A^{p,q}(X)\\
=&\dim H_{BC}^{p,q}(X)-v^{p,q}_t-w^{p,q}_t.
\end{align*}
This is \eqref{eq-BCjump}. For \eqref{eq-Ajump}, since
\[
H_{BC}^{p,q}(X)\cong H_{A}^{n-q,n-p}(X),\quad H_{A}^{p,q}(X)\cong H_{A}^{q,p}(X),\quad \forall (p,q)\in \mathbb{N}\times \mathbb{N}
\]
and $H_{BC\phi(t)}^{p,q}(X)\cong H_{A\phi(t)}^{n-p,n-q}(X)$ by Corollary \ref{coro-BCphi-Aphi}, it follows from \eqref{eq-BCjump} that for each $(p,q)\in \mathbb{N}\times \mathbb{N}$,
\begin{align*}
\dim H_{A}^{n-p,n-q}(X)=&\dim H_{A}^{n-q,n-p}(X)=\dim H_{BC}^{p,q}(X)\\
=&\dim H_{BC\phi(t)}^{p,q}(X)+v^{p,q}_t+w^{p-1,q-1}_t\\
=&\dim H_{A\phi(t)}^{n-p,n-q}(X)+v^{p,q}_t+w^{p-1,q-1}_t.
\end{align*}
\eqref{eq-Ajump} follows from this by replacing $p, q$ with $n-p, n-q$ respectively.
\end{proof}
Note that since the Aeppli deformation of $(n,q)$ or $(p,n)$-forms is always canonically unobstructed, we have $w^{n,q}_t=w^{p,n}_t=0$ for any $t\in B$.

The jumping formulas provide us a useful criterion to determine when the dimension of these deformed cohomology remains constant under deformation of complex structures. In fact, set
\begin{align*}
h_{BC\phi(t)}^{p,q}=&\dim H_{BC\phi(t)}^{p,q}(X),\quad h_{BC}^{p,q}=\dim H_{BC}^{p,q}(X),\\
h_{A\phi(t)}^{p,q}=&\dim H_{A\phi(t)}^{p,q}(X),\quad h_{A}^{p,q}=\dim H_{A}^{p,q}(X),
\end{align*}
we have the following
\begin{corollary}In the situation of Theorem \ref{thm-BCandA-jump-formula}.
\begin{enumerate}
  \item The deformed Bott-Chern number $h_{BC\phi(t)}^{p,q}$ is independent of $t\in B$ if and only if the Bott-Chern deformations of $(p,q)$-forms and the Aeppli deformations of $(p-1,q-1)$-forms are canonically unobstructed;
  \item The deformed Aeppli number $h_{A\phi(t)}^{p,q}$ is independent of $t\in B$ if and only if the Bott-Chern deformations of $(n-p,n-q)$-forms and the Aeppli deformations of $(n-p-1,n-q-1)$-forms are canonically unobstructed;
\end{enumerate}
\end{corollary}
\begin{proof}This follows directly from Definition \eqref{def-unobstructedness}, the jumping formula \eqref{eq-BCjump} and \eqref{eq-Ajump}.
\end{proof}

\subsection{Semi-continuous property under the analytic Zariski topology}
\begin{proposition}\label{Prop-upper-continuous-vw}
For each $(p,q)\in \mathbb{N}\times \mathbb{N}$, $v^{p,q}_t$ and $w^{p,q}_t$ defined in Theorem \ref{thm-BCandA-jump-formula} are lower semi-continuous functions of $t$ under the analytic Zariski topology.
\end{proposition}
\begin{proof}The case for $v^{p,q}_t$ is already known, see \cite[Remark\,4.9]{Xia19dBC}. So we only need prove this for $w^{p,q}_t$. Indeed, since $W_{t}^{p,q}$ is isomorphic to $\ker d^*\cap\ker \p\pb_{\phi(t)} \cap A^{p,q}(X)$, we have
\[
w^{p,q}_t=\dim H_{A}^{p,q}(X)-\dim W_{t}^{p,q}.
\]
On the other hand, assume $\{e_1,\cdots, e_m\}$ is an orthonormal basis of $\mathcal{H}_{A}^{p+1,q+1}(X)$ and $V=\mathcal{H}_{A}^{p,q}(X)$. Then it follows from Proposition \ref{prop-Aclosed=BCharmonic part vanish} that for any $k\in\N$
\begin{align*}
W_{t}^{p,q}=
&\{\sigma_0\in \mathcal{H}_{A}^{p,q}(X) \mid  \p\pb_{\phi(t)}\sigma(t)=0\}\\
=&\{\sigma_0\in \mathcal{H}_{A}^{p,q}(X) \mid  \mathcal{H}_{BC}\Big(\p i_{\phi(t)}\p\sigma(t)\Big)=0\}\\
=&\{\sigma_0\in \mathcal{H}_{A}^{p,q}(X) \mid  \langle \p i_{\phi(t)}\p\sigma(t), e_i\rangle_{L^2}=0,~1\leq i\leq m\}\\
=&\Big\{\sum_{l=1}^N a_l\sigma_0^l\in \mathcal{H}_{A}^{p,q}(X) \mid  \sum_{l=1}^N a_l\langle \p i_{\phi(t)}\p\sigma^l(t), e_i\rangle_{L^2}=0,~1\leq i\leq m\Big\}
\end{align*}
where $\{\sigma_0^l\}_{l=1}^N$ is a basis of $\mathcal{H}_{A}^{p,q}(X)$ and $\sigma (t)=\sum_{k} \sigma_k$ is the canonical deformation of $\sigma_0$ with $\sigma_k=G_{A}(\p\bar{\p})^*\sum_{i+j=k} \p i_{\phi_j} \p\sigma_i,~\forall k>0$. Similarly, each $\sigma^l(t)$ is the canonical deformation of $\sigma_0^l$. As a result, we get
\[
\{t\in B\mid w^{p,q}_t< k\}
=\Big\{t\in B\mid\dim H_{A}^{p,q}(X)-\dim W_{t}^{p,q}=\textrm{rank}\Big(u_i^l(t)\Big)_{m\times N}< k\Big\},
\]
where $u_i^l(t)=\langle \p i_{\phi(t)}\p\sigma^l(t), e_i\rangle_{L^2}$. Because each $u_i^l(t)$ is holomorphic in $t$, we see that for any $k\in\N$, the set $\{t\in B\mid w^{p,q}_t< k\}$ is an analytic subset of $B$.
\end{proof}

\begin{proposition}\label{Prop-upper-continuous-alt-sum}In the situation of Theorem \ref{thm-BCandA-jump-formula}. For any $(p,q)\in \mathbb{N}\times \mathbb{N}$, the following alternating sums are upper semi-continuous functions of $t$ under the analytic Zariski topology:
  \begin{enumerate}
  \item $\sum_{i=0}^q(-1)^{q-i}(h_{BC\phi(t)}^{p,i}+v^{p,i}_t-w^{p-1,i}_t)$;
  \item $\sum_{i=0}^q(-1)^{q-i}(h_{A\phi(t)}^{p,i}+v^{n-p,n-i}_t-w^{n-p-1,n-i}_t)$.
  \end{enumerate}
\end{proposition}
\begin{proof}We will only give the proof of $(1)$ since $(2)$ can be shown by the exactly the same method. For $(1)$, by applying the jumping formula \eqref{eq-BCjump}, we get
\begin{align*}
\begin{split}
  h_{BC}^{p,q}=&(h_{BC\phi(t)}^{p,q}+v^{p,q}_t-w^{p-1,q}_t)+w^{p-1,q}_t+w^{p-1,q-1}_t,\\
  -h_{BC}^{p,q-1}=&-(h_{BC\phi(t)}^{p,q-1}+v^{p,q-1}_t-w^{p-1,q-1}_t)-w^{p-1,q-1}_t-w^{p-1,q-2}_t,\\
   &\vdots\\
  (-1)^q h_{BC}^{p,0}=&(-1)^q(h_{BC\phi(t)}^{p,0}+v^{p,0}_t-w^{p-1,0}_t)+(-1)^qw^{p-1,0}_t.\\
 \end{split}
\end{align*}
Summing all these up, we have
\[
\sum_{i=0}^q(-1)^{q-i} h_{BC}^{p,i}-w^{p-1,q}_t=\sum_{i=0}^q(-1)^{q-i}(h_{BC\phi(t)}^{p,i}+v^{p,i}_t-w^{p-1,i}_t).
\]
The conclusion then follows since by Proposition \ref{Prop-upper-continuous-vw} we know $w^{p-1,q}_t$ is a lower semi-continuous functions of $t$ under the analytic Zariski topology.
\end{proof}
\begin{remark}
 Similarly, we can show $\sum_{i=0}^q(-1)^{q-i}(h_{BC\phi(t)}^{p,i}-v^{p,i-1}_t+w^{p-1,i-1}_t)$ and
$\sum_{i=0}^q(-1)^{q-i}(h_{A\phi(t)}^{p,i}-v^{n-p,n-i+1}_t+w^{n-p-1,n-i-1}_t)$ is upper semi-continuous.
\end{remark}
\subsection{Remarks on canonical Dolbeault deformations}
We end this section with some discussions about canonical Dolbeault deformations.

Let $E$ be a holomorphic tensor bundle on $X$ and the $E_t$ the corresponding tensor bundle on $X$. It is proved \cite[Thm.\,5.10]{Xia19dDol} by the second author that
\begin{equation}\label{eq-vtpq}
\dim H^{q}(X,E)=\dim H^{q}(X_t,E_t)+v^{q}_t+v^{q-1}_t,
\end{equation}
where $v^{q}_t:=\dim H^{q}(X,E)-\dim \ker\pb_{\phi(t)}\cap\ker\pb^*\cap A^{0,q}(X,E) \geq 0$. By using similar arguments as in the proof of Proposition \ref{Prop-upper-continuous-alt-sum}, it follows from \eqref{eq-vtpq} that
\begin{equation}\label{eq-vtq=altsumhi}
v^{q}_t=\sum_{i=0}^{q}(-1)^{q-i}\dim H^{i}(X,E)-\sum_{i=0}^{q}(-1)^{q-i}\dim H^{i}(X_t,E_t).
\end{equation}
Since the Dolbeault deformations of $E$-valued $(0,q)$-forms are canonically unobstructed iff $v^{q}_t=0$ for any $t\in B$, we conclude that for Dolbeault deformations of $E$-valued $(0,q)$-forms,
\[
\text{canonically~unobstructedness}\Longleftrightarrow\sum_{i=0}^{q}(-1)^{q-i}\dim H^{i}(X_t,E_t)~\text{is~independent~of}~t.
\]
In particular, we get the following
\begin{corollary}Given a Hermitian metric $h$ on $X$, assume the deformations of $E$-valued $(0,q)$-forms is canonically unobstructed with respect to $h$. Then the deformations of $E$-valued $(0,q)$-forms is also canonically unobstructed with respect to any other Hermitian metrics on $X$.
\end{corollary}
As a result, the notion of canonically unobstructedness (at least for Dolbeault deformations) is independent of the choices of Hermitian metrics.

On the other hand, we see that $\sum_{i=0}^{q}(-1)^{q-i}\dim H^{i}(X_t,E_t)$ is upper semicontinuous (under the analytic Zariski topology) since $v^{q}_t$ is lower semi-continuous. This is in fact a special case of a result due to Flenner \cite{Fle81} where $\{E_t\}_{t\in B}$ may be replaced by a coherent sheaf on the total space $\mathcal{X}$, see also \cite[Thm.\,5.10]{BDIP02}.

\section{Unobstructed Bott-Chern/Aeppli deformations}
In this section, we study when the Bott-Chern/Aeppli deformations of $(p,q)$-forms are (canonically) unobstructed.

Let $\pi: (\mathcal{X}, X)\to (B,0)$ be a complex analytic family with Beltrami differentials $\phi(t)$. Write $\phi(t)=\sum_j\phi_j$, then the Maurer-Cartan equation $\pb\phi(t)=\frac{1}{2}[\phi(t),\phi(t)],~t\in B$ implies
\begin{equation}\label{eq-Maurer-Cartan-expanded}
\pb\phi_k=\frac{1}{2}\sum_{j=1}^k[\phi_j,\phi_{k-j}],\quad k\geq 1.
\end{equation}
It follows from this and Cartan formulas (c.f. \cite{LRY15,Xia19deri,Xia19dDol,BM18,FM09,FM06,Cle05}) that
\begin{align*}
i_{\phi_j}\bar{\partial}
=&\bar{\partial}i_{\phi_j}-i_{\bar{\partial}\phi_j}\\
=&\bar{\partial}i_{\phi_j}-\frac{1}{2}\sum_{l=1}^{j}(\mathcal{L}_{\phi_l}^{1,0}i_{\phi_{j-l}}-i_{\phi_{j-l}}\mathcal{L}_{\phi_l}^{1,0})\\
=&\bar{\partial}i_{\phi_j}-\frac{1}{2}\sum_{l=1}^{j}(i_{\phi_l}\partial i_{\phi_{j-l}}-\partial i_{\phi_l}i_{\phi_{j-l}}-i_{\phi_{j-l}}i_{\phi_l}\partial +i_{\phi_{j-l}}\partial i_{\phi_l})
\end{align*}
which implies
\begin{equation}\label{eq-iphi-pb}
i_{\phi_j}\pb=\pb i_{\phi_j}-\sum_{l=1}^{j}i_{\phi_l}\p i_{\phi_{j-l}}+\frac{1}{2}\sum_{l=1}^{j}(\p i_{\phi_l}i_{\phi_{j-l}}+i_{\phi_{j-l}}i_{\phi_l}\p).
\end{equation}

\begin{theorem}\label{thm-unobstructed-deformations}
Let $\pi: (\mathcal{X}, X)\to (B,0)$ be a complex analytic family with Beltrami differentials $\phi(t)$. Suppose $X$ is equipped with a fixed Hermitian metric.
\begin{enumerate}
  \item Assume $\p_{A,\pb(\ker\p)}^{p-1,q+1}=0$. Then the Bott-Chern deformations of $(p,q)$-forms are canonically unobstructed;
  \item Assume $\p_{A,BC}^{p,q+1}=0$. Then the Aeppli deformations of $(p,q)$-forms are canonically unobstructed;
  \item Assume $\p_{A,\pb(\ker\p)}^{p-1,q+1}=0$ and $\p_{A,BC}^{p-1,q}=0$. Then the deformed Bott-Chern number $h_{BC\phi(t)}^{p,q}$ is independent of $t\in B$;
  \item Assume $\p_{A,\pb(\ker\p)}^{n-p-1,n-q+1}=0$ and $\p_{A,BC}^{n-p-1,n-q}=0$. Then the deformed Aeppli number $h_{A\phi(t)}^{p,q}$ is independent of $t\in B$.
\end{enumerate}
\end{theorem}
\begin{proof}$(1)$ It is enough to show the following system of equations
\begin{equation}\label{eq-ob-BC}
\left\{
\begin{array}{ll}
\pb\sigma_k= \sum_{j=1}^{k}\mathcal{L}_{\phi_j}^{1,0}\sigma_{k-j} , &\quad k\geq 1  \\
\p\sigma_k=0, &\quad k\geq 1\\
\end{array} \right.
\end{equation}
can be solved inductively for any given $\sigma_0\in\ker d\cap A^{p,q}(X)$. In fact, assume $\sigma_1, \cdots, \sigma_{N}\in A^{p,q}(X)$ are solutions of \eqref{eq-ob-BC} for $k\leq N$. We need to solve
\[
\pb\sigma_{N+1}= \sum_{j=1}^{N+1}\mathcal{L}_{\phi_j}^{1,0}\sigma_{N+1-j}=-\sum_{j=1}^{N+1}\p i_{\phi_j}\sigma_{N+1-j}
\]
for $\sigma_{N+1}\in\ker\p\cap A^{p,q}(X)$. For this purpose, we observe that by using \eqref{eq-iphi-pb} and the induction assumption, there holds
\begin{align*}
&\pb(\sum_j i_{\phi_j}\sigma_{N+1-j} )\\
=&\sum_{j,l} i_{\phi_j}\mathcal{L}_{\phi_l}^{1,0}\sigma_{N+1-j-l}+\sum_{j,l}i_{\phi_l}\p i_{\phi_{j-l}}\sigma_{N+1-j}-\frac{1}{2}\sum_{j,l}(\p i_{\phi_l}i_{\phi_{j-l}}+i_{\phi_{j-l}}i_{\phi_l}\p)\sigma_{N+1-j}\\
=&-\frac{1}{2}\sum_{j,l}\p i_{\phi_l}i_{\phi_{j-l}}\sigma_{N+1-j}+ \frac{1}{2}\sum_{j,l} i_{\phi_j}i_{\phi_{l}}\p\sigma_{N+1-j-l}\\
=&-\frac{1}{2}\sum_{j,l}\p i_{\phi_l}i_{\phi_{j-l}}\sigma_{N+1-j}.
\end{align*}
As a result, $\p\pb(\sum_j i_{\phi_j}\sigma_{N+1-j})=0$ which by using $\p_{A,\pb(\ker\p)}^{p-1,q+1}=0$ implies there exists $\sigma_{N+1}\in\ker\p\cap A^{p,q}(X)$ such that
\[
\pb\sigma_{N+1}= -\sum_{j=1}^{N+1}\p i_{\phi_j}\sigma_{N+1-j}.
\]

$(2)$ Similarly, we want to show the following system of equations
\begin{equation}\label{eq-ob-A}
\p\pb\sigma_k= \p\sum_{j=1}^{k}\mathcal{L}_{\phi_j}^{1,0}\sigma_{k-j} , \quad k\geq 1,
\end{equation}
can be solved inductively for any given $\sigma_0\in\ker \p\pb\cap A^{p,q}(X)$. Assume $\sigma_1, \cdots, \sigma_{N}\in A^{p,q}(X)$ are solutions of \eqref{eq-ob-A} for $k\leq N$. We need to solve
\[
\p\pb\sigma_{N+1}= \p\sum_{j=1}^{N+1}\mathcal{L}_{\phi_j}^{1,0}\sigma_{N+1-j}
\]
for $\sigma_{N+1}\in A^{p,q}(X)$. In view of the assumption $\p_{A,BC}^{p,q+1}=0$, this can be done if $\sum_{j=1}^{k}\mathcal{L}_{\phi_j}^{1,0}\sigma_{k-j}\in\ker\p\pb$ .

Now we show $\sum_{j=1}^{k}\mathcal{L}_{\phi_j}^{1,0}\sigma_{k-j}\in\ker\p\pb$. In fact, it follows from \eqref{eq-Maurer-Cartan-expanded}, the induction assumption and Cartan formulas that (c.f. the proof of \cite[Thm.\,3.3]{Xia19dDol})
\begin{align*}
&\pb\sum_{j=1}^{N}\mathcal{L}_{\phi_j}^{1,0}\sigma_{N-j}\\
=&\sum_{j=1}^{N} \left( -\mathcal{L}_{\phi_j}^{1,0} \pb\sigma_{N-j} + \mathcal{L}_{\pb\phi_j}^{1,0}\sigma_{N-j}  \right)\\
=&\sum_{j=1}^{N} \left( -\mathcal{L}_{\phi_j}^{1,0} \sum_{i=1}^{N-j}\mathcal{L}_{\phi_i}^{1,0} \sigma_{N-j-i} +
\frac{1}{2}\mathcal{L}_{ \sum_{i=1}^{j}[\phi_i, \phi_{j-i}]}^{1,0} \sigma_{N-j} \right) -\sum_{j=1}^{N}\sum_{i=1}^{N-j}\mathcal{L}_{\phi_j}^{1,0} \left(\pb-\mathcal{L}_{\phi_i}^{1,0}\right) \sigma_{N-j-i}\\
=&-\sum_{j=1}^{N}\sum_{i=1}^{N-j}\mathcal{L}_{\phi_j}^{1,0} \left(\pb-\mathcal{L}_{\phi_i}^{1,0}\right) \sigma_{N-j-i}
=\sum_{j=1}^{N}\sum_{i=1}^{N-j}\p i_{\phi_j} \left(\pb-\mathcal{L}_{\phi_i}^{1,0}\right) \sigma_{N-j-i}\in\ker\p.
\end{align*}
Finally, $(3)$ follows from $(1)$ and \eqref{eq-BCjump}. $(4)$ is restatement of $(3)$ in view of \eqref{eq-BC-A-dual}.
\end{proof}

\begin{corollary}
The Bott-Chern number $\dim H_{BC}^{p,0}(X_t)=\dim H_{A}^{n,n-p}(X_t)$ is independent of $t\in B$ if $\p_{A,\pb(\ker\p)}^{p-1,1}=0$.
\end{corollary}
\begin{proof}Since $w_t^{p-1,-1}=0$ for any $t\in B$. When $q=0$, the equality \eqref{eq-BCjump} becomes
\[
\dim H_{BC}^{p,0}(X)=\dim H_{BC\phi(t)}^{p,0}(X)+v^{p,0}_t.
\]
The conclusion now follows from this and $(1)$ of Theorem \ref{thm-unobstructed-deformations}.
\end{proof}

\section{Two concrete examples}\label{sec-examples}
In this section, we present explicit computations of canonical Bott-Chern/Aeppli deformations to examine our results. We will consider the Iwasawa manifold and the Nakamura manifold. Both of them are holomorphically parallelizable, i.e. of the form $\Gamma\setminus G$ where $G$ is a complex Lie group and $\Gamma\subset G$ is a cocompact lattice of maximal rank. For the Iwasawa manifold, $G$ is nilpotent. So according to \cite{Sak76}, combined with \cite{Ang13} or \cite[Coro.\,13]{Ste21}, the Dolbeault/Bott-Chern/Aeppli cohomology can be computed by (left) invariant forms. For the Nakamura manifold, $G$ is not nilpotent but only solvable. In this case, we need to apply results of Angella-Kasuya \cite{AK17a,AK17b} to compute the Bott-Chern/Aeppli cohomology.
\subsection{The Iwasawa manifold}
\begin{example}\label{example Case III-(2)}
Case III-(2). Let $G$ be the matrix Lie group defined by
\[
G := \left\{
\left(
\begin{array}{ccc}
 1 & z^1 & z^3 \\
 0 &  1  & z^2 \\
 0 &  0  &  1
\end{array}
\right) \in \mathrm{GL}(3;\mathbb{C}) \mid z^1,\,z^2,\,z^3 \in\mathbb{C} \right\}\cong \mathbb{C}^3,
\]
where the product is the one induced by matrix multiplication. This is usually called the \emph{Heisenberg group}. Consider the discrete subgroup $\Gamma$ defined by
\[
\Gamma := \left\{
\left(
\begin{array}{ccc}
 1 & \omega^1 & \omega^3 \\
 0 &  1  & \omega^2 \\
 0 &  0  &  1
\end{array}
\right) \in G \mid \omega^1,\,\omega^2,\,\omega^3 \in\mathbb{Z}[\sqrt{-1}] \right\},
\]
The quotient $X=\Gamma\setminus G$ is called the \emph{Iwasawa manifold}. A basis of $H^0(X,\Omega^1)$ is given by
\[
\varphi^1 = d z^1,~ \varphi^2 = d z^2,~ \varphi^3 = d z^3-z^1\,d z^2,
\]
and a dual basis $\theta^1, \theta^2, \theta^3\in H^0(X,T_X^{1,0})$ is given by
\[
\theta^1=\frac{\partial}{\partial z^1},~\theta^2=\frac{\partial}{\partial z^2} + z^1\frac{\partial}{\partial z^3},~\theta^3=\frac{\partial}{\partial z^3}.
\]
The structure equation in terms of the basis $\{\varphi^1,\varphi^2,\varphi^3\}$ is given by
\[
d\varphi^1=0,\quad d\varphi^2=0,\quad d\varphi^3=-\varphi^{12}.
\]
$X$ is equipped with the Hermitian metric $\sum_{i=1}^3\varphi^i\otimes\bar{\varphi}^i$. The Beltrami differential of the Kuranishi family of $X$ is
\[
\phi(t) = \sum_{i=1}^3\sum_{\lambda=1}^2t_{i\lambda}\theta^i\bar{\varphi}^{\lambda} - D(t)\theta^3\bar{\varphi}^{3},~\text{with}~D(t)=t_{11}t_{22}-t_{21}t_{12},
\]
and the Kuranishi space of $X$ is
\[
\mathcal{B}=\{t=(t_{11}, t_{12}, t_{21}, t_{22}, t_{31}, t_{32})\in \mathbb{C}^6\mid |t_{i\lambda}|<\epsilon, i=1, 2, 3, \lambda=1,2 \},
\]
where $\epsilon>0$ is sufficiently small. Set
\[
\phi_1=\sum_{i=1}^3\sum_{\lambda=1}^2t_{i\lambda}\theta^i\bar{\varphi}^{\lambda},~ \phi_2 = -D(t)\theta^3\bar{\varphi}^{3},
\]
then by direct computations, we have
\begin{align*}
&\mathcal{L}_{\phi_1}^{1,0}\varphi^1 = \mathcal{L}_{\phi_1}^{1,0}\varphi^2=\mathcal{L}_{\phi_1}^{1,0}\bar{\varphi}^1=\mathcal{L}_{\phi_1}^{1,0}\bar{\varphi}^{2}=\mathcal{L}_{\phi_1}^{1,0}\bar{\varphi}^{3}=0,\\[5pt]
&\mathcal{L}_{\phi_1}^{1,0}\varphi^3 = \sum_{\lambda=1}^2(t_{1\lambda}\varphi^2-t_{2\lambda}\varphi^1)\wedge\bar{\varphi}^{\lambda},~\mathcal{L}_{\phi_2}^{1,0}\varphi^i =\mathcal{L}_{\phi_2}^{1,0}\bar{\varphi}^{i}=0,~i=1,2,3.
\end{align*}
In particular,
\begin{equation}\label{eq-plie}
\p\mathcal{L}_{\phi_1}^{1,0}\varphi^i=\p\mathcal{L}_{\phi_1}^{1,0}\bar{\varphi}^i=\mathcal{L}_{\phi_2}^{1,0}\varphi^i =\mathcal{L}_{\phi_2}^{1,0}\bar{\varphi}^{i}=0,\quad i=1, 2, 3.
\end{equation}
\end{example}
The Bott-Chern/Aeppli cohomology of $X$ are easily computed by using left invariant forms. Let us first consider canonical Aeppli deformations.
\begin{itemize}
    \item
   $H_{A}^{1,0}(X)=\mathbb{C}\{\varphi^1,\varphi^2,\varphi^3\}$,\quad $H_{A}^{0,1}(X)=\mathbb{C}\{\varphi^{\bar{1}},\varphi^{\bar{2}},\varphi^{\bar{3}}\}$;
   \item
   $H_{A}^{2,0}(X)=\mathbb{C}\{\varphi^{13},\varphi^{23} \}$,\quad  $H_{A}^{0,2}(X)=\mathbb{C}\{\varphi^{\overline{13}},\varphi^{\overline{23}}\}$;
   \item
  $H_{A}^{1,1}(X)=\mathbb{C}\{\varphi^{1\bar{1}},\varphi^{1\bar{2}},\varphi^{1\bar{3}},\varphi^{2\bar{1}},\varphi^{2\bar{2}},\varphi^{2\bar{3}},\varphi^{3\bar{1}},\varphi^{3\bar{2}}\}$;
   \item
   $H_{A}^{3,0}(X)=\mathbb{C}\{\varphi^{123}\}$,\quad$H_{A}^{2,1}(X)=\mathbb{C}\{\varphi^{13\bar{1}},\varphi^{13\bar{2}},\varphi^{13\bar{3}},\varphi^{23\bar{1}},\varphi^{23\bar{2}},\varphi^{23\bar{3}}\}$;
   \item
   $H_{A}^{1,2}(X)=\mathbb{C}\{\varphi^{1\overline{13}},\varphi^{1\overline{23}},\varphi^{2\overline{13}},\varphi^{2\overline{23}},\varphi^{3\overline{13}},\varphi^{3\overline{23}}\}$,\quad$H_{A}^{0,3}(X)=\mathbb{C}\{\varphi^{\overline{123}}\}$;
   \item
   $H_{A}^{3,1}(X)=\mathbb{C}\{\varphi^{123\bar{1}},\varphi^{123\bar{2}},\varphi^{123\bar{3}}\}$,\quad$H_{A}^{2,2}(X)=\mathbb{C}\{\varphi^{13\overline{23}},\varphi^{23\overline{13}},\varphi^{23\overline{23}},\varphi^{13\overline{13}}\}$
   \item
   $H_{A}^{1,3}(X)=\mathbb{C}\{\varphi^{1\overline{123}},\varphi^{2\overline{123}},\varphi^{3\overline{123}}\}$
   \item
   $H_{A}^{3,2}(X)=\mathbb{C}\{\varphi^{123\overline{13}},\varphi^{123\overline{23}}\}$,\quad$H_{A}^{2,3}(X)=\mathbb{C}\{\varphi^{13\overline{123}},\varphi^{23\overline{123}}\}$
   \item
   $H_{A}^{3,3}(X)=\mathbb{C}\{\varphi^{123\overline{123}}\}$
\end{itemize}
In fact, in view of \eqref{eq-plie}, for any $\sigma_0\in\ker\p\pb\cap A^{p,q}(X)$, we have
\[
\p i_{\phi_1}\p\sigma_0=\p\mathcal{L}_{\phi_1}^{1,0}\sigma_0=0
\]
and thus
\[
\sigma_1 = G_{A}(\p\bar{\p})^* \p i_{\phi_1}\p \sigma_0=0,\quad  \sigma_2 = G_{A}(\p\bar{\p})^* (\p i_{\phi_2}\p \sigma_0 + \p i_{\phi_1}\p \sigma_1 )=0.
\]
Similarly, there holds
\[
\p\sum_{j=1}^{k}\mathcal{L}_{\phi_j}^{1,0}\sigma_{k-j}=0 , \quad k\geq 1,
\]
and $\sigma_k= G_{A}(\p\bar{\p})^*\sum_{i+j=k} \p i_{\phi_j} \p\sigma_i=0$ for any $k\geq1$. Therefore, for each $(p,q)\in \mathbb{N}\times \mathbb{N}$ the Aeppli deformations of $(p,q)$-forms are canonically unobstructed and $w_t^{p,q}=0$ for any $t\in\mathcal{B}$.

Let us now consider canonical Bott-Chern deformations.
\begin{itemize}
   \item $H_{BC}^{2,0}(X)=\mathbb{C}\{\varphi^{12},\varphi^{13},\varphi^{23} \}$;
   \item $H_{BC}^{2,1}(X)=\mathbb{C}\{\varphi^{12\bar{1}},\varphi^{12\bar{2}},\varphi^{13\bar{1}},\varphi^{13\bar{2}},\varphi^{23\bar{1}},\varphi^{23\bar{2}}\}$;
   \item $H_{BC}^{2,2}(X)=\mathbb{C}\{\varphi^{12\overline{23}},\varphi^{23\overline{23}},\varphi^{13\overline{23}},
\varphi^{12\overline{13}},\varphi^{23\overline{13}},\varphi^{13\overline{13}},\varphi^{23\overline{12}},\varphi^{13\overline{12}}\}$.
\end{itemize}
For $(p,q)=(2,0)$. Let $\sigma_0 = a_{12}\varphi^{12} +a_{13}\varphi^{13}+a_{23}\varphi^{23}$, we have
\[
\p i_{\phi_1}\sigma_0=(a_{13}t_{11}+a_{23}t_{21})\varphi^{12\bar{1}}+(a_{13}t_{12}+a_{23}t_{22})\varphi^{12\bar{2}}
\]
is exact if and only if $\p i_{\phi_1}\sigma_0=0$, i.e.
\begin{equation}\label{a_{23}}
\left\{
\begin{array}{rcl}
t_{21}a_{23} + t_{11}a_{13} &=& 0 \\[5pt]
t_{22}a_{23} + t_{12}a_{13} &=& 0.
\end{array}
\right.
\end{equation}
has solutions for $(a_{23}, a_{13})$, and in this case the (degree $1$ term of) canonical deformation is given by
\[
\sigma_{1}=-G_{BC}(\pb^*\p\p^*+\pb^*)\partial i_{\phi_1}\sigma_0 =0.
\]
On the other hand,
\[
\partial i_{\phi_2}\sigma_0 = 0 \Longrightarrow \sigma_{2}=-G_{BC}(\pb^*\p\p^*+\pb^*)\partial(i_{\phi_2}\sigma_0+i_{\phi_1}\sigma_1) =0,
\]
and $\phi_k=0,~ k>2$  implies that $\sigma_{k}= G_{BC}(\pb^*\p\p^*+\pb^*)\sum_{i+j=k} \p i_{\phi_j}\sigma_i=0$ for any $k>2$.

Therefore, for $V=H_{BC}^{2,0}(X)$ we have (See \eqref{eq-V BC,t})
\[
V_{BC,t}^{2,0}=\{a_{12}\varphi^{12} + a_{23}\varphi^{23} + a_{13}\varphi^{13} \mid (a_{12}, a_{23}, a_{13})\in \mathbb{C}^3~\text{satisfy}~ \eqref{a_{23}}\}\subseteq H_{BC}^{2,0}(X),
\]
and
\[
v_t^{2,0}=\dim H_{BC}^{2,0}(X)-\dim V_{BC,t}^{2,0}=
\left\{
\begin{array}{rcl}
0, &(t_{11},t_{12},t_{21},t_{22}) =0 \\[5pt]
1, &(t_{11},t_{12},t_{21},t_{22}) \neq 0,D(t)=0 \\[5pt]
2, &(t_{11},t_{12},t_{21},t_{22}) \neq 0,D(t)\neq 0.
\end{array}
\right.
\]
For $(p,q)$=(2,1). Let $\sigma_0 = \sum a_{ij{\bar{k}}}\varphi^{ij{\bar{k}}}\in\mathcal{H}_{BC}^{2,1} $, we have
\begin{align*}
\p i_{\phi_1}\sigma_0&=(t_{11}a_{13\bar{2}}-t_{12}a_{13\bar{1}}+t_{21}a_{23\bar{2}}-t_{22}a_{23\bar{1}})\varphi^{12\overline{12}}\\
&=a(t)\p\pb\varphi^{3\bar{3}},
\end{align*}
where $a(t):=-(t_{11}a_{13\bar{2}}-t_{12}a_{13\bar{1}}+t_{21}a_{23\bar{2}}-t_{22}a_{23\bar{1}})$. In this case,
\begin{align*}
\sigma_{1}&=-G_{BC}(\pb^*\p\p^*+\pb^*)\p i_{\phi_1}\sigma_0\\
&=-a(t)G_{BC}(\pb^*\p\p^*+\pb^*)\p\pb\varphi^{3\bar{3}}\\
&=a(t)G_{BC}(\pb^*\p\p^*\pb+\pb^*\pb)\p\varphi^{3\bar{3}}\\
&=a(t)G_{BC}(\square_{BC}-(\p\pb)(\p\pb)^*)\p\varphi^{3\bar{3}}\\
&=-a(t)G_{BC}(\square_{BC}-(\p\pb)(\p\pb)^*)\varphi^{12\bar{3}}\\
&=-a(t)G_{BC}\square_{BC}\varphi^{12\bar{3}}\\
&=-a(t)(1-\mathcal{H}_{BC}^{2,1})\varphi^{12\bar{3}}\\
&=-a(t)\varphi^{12\bar{3}}.
\end{align*}
On the other hand, $\p i_{\phi_2}\sigma_0=-\mathcal{L}_{\phi_2}^{1,0}\sigma_0=0$ and $\p i_{\phi_1}\sigma_1=-\mathcal{L}_{\phi_1}^{1,0}\sigma_1=0$. Therefore,
\[
\sigma_{2}=-G_{BC}(\pb^*\p\p^*+\pb^*)\p(i_{\phi_2}\sigma_0+i_{\phi_1}\sigma_1) =0,
\]
and $\phi_k=0,~ k>2$  implies that $\sigma_{k}=0,~ k>2$. It follows that the Bott-Chern deformations of $(2,1)$-forms are canonically unobstructed and $v_t^{2,1}=0$ for any $t\in\mathcal{B}$.

Similarly, we can determine $v_t^{p,q}$ for all $(p,q)$. It turns out there are only two cases for which $v_t^{p,q}\neq0$: $(p,q)=(2,0), (2,2)$ (see \cite[Sec.\,6]{Xia19dBC} for the computations of $v_t^{2,2}$).  As a result, we get the following table (write $(i), (ii), (iii)$ for the three cases when $(t_{11}, t_{12}, t_{21}, t_{22})=0$, $(t_{11}, t_{12}, t_{21}, t_{22})\neq0$ and $D(t)= 0$, $D(t)\neq 0$, respectively):
\renewcommand\arraystretch{1.5}
\begin{table}[!htbp]
\caption{$v^{\bullet,\bullet}_t$ and $w^{\bullet,\bullet}_t$ for cases $(i), (ii), (iii)$}
\centering
\begin{center}
\begin{tabular}{|c|c|c|c|c|c|c|c|}
\hline
  & $(i)$ & $(ii)$ & $(iii)$ & & $(i)$ & $(ii)$ & $(iii)$  \\
\hline
 $h^{1,3}_A$ & $3$ & $3$ & $3$ &$h^{1,1}_A$& $8$ & $8$ & $8$  \\
\hline
 $h^{2,0}_{BC}$ & $3$ & $3$ & $3$ &$h^{2,2}_{BC}$& $8$ & $8$ & $8$  \\
\hline
 $w^{1,-1}_t$ & $0$ & $0$ & $0$ & $w^{1,1}_t$ & $0$ & $0$ & $0$  \\
\hline
 $v^{2,0}_t$ & $0$ & $1$ & $2$ & $v^{2,2}_t$& $0$ & $1$ & $1$  \\
\hline
 $h^{2,0}_{BC,\phi(t)}$ & $3$ & $2$ & $1$ &$h^{2,2}_{BC,\phi(t)}$ & $8$ & $7$ & $7$  \\
\hline
 $h^{1,3}_{A,\phi(t)}$ & $3$ & $2$ & $1$ &$h^{1,1}_{A,\phi(t)}$ & $8$ & $7$ & $7$  \\
\hline
\end{tabular}
\end{center}
\end{table}

\subsection{The complex parallelizable Nakamura manifold}
\begin{example}
Let $X=\Gamma\setminus G$ be the solvable manifold constructed by Nakamura in Example III-(3b) of~\cite{Nak75}. We will use an equivalent description of $X$ given by\footnote{We learned of this construction through communications with D. Angella and H. Kasuya.} De Bartolomeis-Tomassini \cite{DT06} (see also \cite{AK17a,AK17b,Has10}). In fact, following the conventions in \cite[Sec.\,4]{AK17b}, let $G:=(\C^3,*)$ be the complex solvable Lie group with the multiplication $*$ defined by
\[
(z_1,z_2,z_3)*(w_1,w_2,w_3)=(z_1+w_1,e^{w_1}z_2+w_2,e^{-w_1}z_3+w_3).
\]
In other words, $G$ is the semi-direct product $\C\ltimes_\varphi\C^2$ where the homomorphism $\varphi:\C\to\C^2$ is given by
\[
\varphi(z)=\left ( \begin{matrix}
e^{z} & 0\\
0 & e^{-z}\\
\end{matrix} \right ),\quad \forall z\in\C.
\]
The lattice $\Gamma\subset G$ has the following generators
\begin{align*}
T_1(z)&:=(z_1+\lambda,e^{\lambda}z_2,e^{-\lambda}z_3),\quad T_2(z):=(z_1+2\pi i,z_2,z_3),\\
T_3(z)&:=(z_1,z_2+1,z_3-\mu),\quad T_4(z):=(z_1,z_2+\mu,z_3+1),\\
T_5(z)&:=(z_1,z_2+2\pi i,z_3-2\pi i\mu),\quad T_6(z):=(z_1,z_2+2\pi i\mu,z_3+2\pi i),
\end{align*}
where $\lambda=\log\frac{3+\sqrt{5}}{2}$ and $\mu=\log\frac{\sqrt{5}-1}{2}$ (c.f. \cite{DT06}).
We will consider the deformation family of $X$ whose Beltrami differential is
\[
\phi(t)=\phi_1=t\frac{\p}{\p z^{1}}dz^{1}.
\]
\end{example}

According to Corollary $2.21$ and Theorem $2.22$ in \cite{AK17b}, there is a sub-double complex $C_\Gamma^{\bullet,\bullet}$ (explicitly determined) of $(A^{\bullet,\bullet}(X), d=\p+\pb)$ such that
\begin{equation}\label{eq-Cgamma-BC}
H^{\bullet,\bullet}_{\pb}(C_\Gamma^{\bullet,\bullet})\cong H^{\bullet,\bullet}_{\pb}(X),\quad H^{\bullet,\bullet}_{BC}(C_\Gamma^{\bullet,\bullet})\cong H^{\bullet,\bullet}_{BC}(X).
\end{equation}
Combining this with \cite[Coro.\,13]{Ste21}, we also have
\begin{equation}\label{eq-Cgamma-A}
H^{\bullet,\bullet}_{A}(C_\Gamma^{\bullet,\bullet})\cong H^{\bullet,\bullet}_{A}(X).
\end{equation}
Thanks to the works of Angella-Kasuya, $C_\Gamma^{\bullet,\bullet}$ has already been computed in \cite[Sec.\,4]{AK17b}. In fact, we have
\begin{align*}
C_{\Gamma}^{0,0}=&\C\{1\},\quad C_{\Gamma}^{3,0}=\mathbb{C}\{dz^{123}\},\quad C_{\Gamma}^{0,3}=\mathbb{C}\{dz^{\overline{123}}\},\quad C_{\Gamma}^{3,3}=\mathbb{C}\{dz^{123\overline{123}}\}\\
C_{\Gamma}^{1,0}=&\mathbb{C}\{dz^1,e^{-z^1}dz^2,e^{z^1}dz^3,e^{-z^{\bar{1}}}dz^2,e^{z^{\bar{1}}}dz^3\}\\
C_{\Gamma}^{0,1}=&\mathbb{C}\{dz^{\bar{1}},e^{-z^1}dz^{\bar{2}},e^{z^1}dz^{\bar{3}},e^{-z^{\bar{1}}}dz^{\bar{2}},e^{z^{\bar{1}}}dz^{\bar{3}}\}\\
C_{\Gamma}^{2,0}=&\mathbb{C}\{e^{-z^{1}}dz^{12},e^{z^{1}}dz^{13},dz^{23},e^{-z^{\bar{1}}}dz^{12},e^{z^{\bar{1}}}dz^{13} \}\\
C_{\Gamma}^{1,1}=&\mathbb{C}\{dz^{1\bar{1}},e^{-z^{1}}dz^{1\bar{2}},e^{z^{1}}dz^{1\bar{3}},e^{-z^{1}}dz^{2\bar{1}},e^{-2z^{1}}dz^{2\bar{2}},dz^{2\bar{3}},e^{z^{1}}dz^{3\bar{1}},dz^{3\bar{2}},\\
&e^{2z^{1}}dz^{3\bar{3}},e^{-z^{\bar{1}}}dz^{2\bar{1}},e^{-z^{\bar{1}}}dz^{1\bar{2}},e^{z^{\bar{1}}}dz^{1\bar{3}},e^{z^{1}}dz^{3\bar{1}},e^{-2z^{\bar{1}}}dz^{2\bar{2}},e^{2z^{\bar{1}}}dz^{3\bar{3}} \}\\
C_{\Gamma}^{0,2}=&\mathbb{C}\{e^{-z^{1}}dz^{\overline{12}},e^{z^{1}}dz^{\overline{13}},dz^{\overline{23}},e^{-z^{\bar{1}}}dz^{\overline{12}},e^{z^{\bar{1}}}dz^{\overline{13}}\}\\
C_{\Gamma}^{2,1}=&\mathbb{C}\{e^{-z^{1}}dz^{12\bar{1}},e^{-2z^{1}}dz^{12\bar{2}},dz^{12\bar{3}},e^{z^{1}}dz^{13\bar{1}},dz^{13\bar{2}},e^{2z^{1}}dz^{13\bar{3}},dz^{23\bar{1}},\\
&e^{-z^{1}}dz^{23\bar{2}},e^{z^{1}}dz^{23\bar{3}},e^{-z^{\bar{1}}}dz^{12\bar{1}},e^{z^{\bar{1}}}dz^{13\bar{1}},e^{-2z^{\bar{1}}}dz^{12\bar{2}},e^{-z^{\bar{1}}}dz^{23\bar{2}},e^{2z^{\bar{1}}}dz^{13\bar{3}},e^{z^{\bar{1}}}dz^{23\bar{3}}  \}\\
C_{\Gamma}^{1,2}=&\mathbb{C}\{e^{-z^{\bar{1}}}dz^{1\overline{12}},e^{-2z^{\bar{1}}}dz^{2\overline{12}},dz^{3\overline{12}}, e^{z^{\bar{1}}}dz^{1\overline{13}},dz^{2\overline{13}},e^{2z^{\bar{1}}}dz^{3\overline{13}},dz^{1\overline{23}},\\
&e^{-z^{\bar{1}}}dz^{2\overline{23}},e^{z^{\bar{1}}}dz^{3\overline{23}},e^{-z^{1}}dz^{1\overline{12}},e^{z^{1}}dz^{1\overline{13}},e^{-2z^{1}}dz^{2\overline{12}},e^{-z^{1}}dz^{2\overline{23}},e^{2z^{1}}dz^{3\overline{13}},e^{z^{1}}dz^{3\overline{13}} \}
\end{align*}
and
\begin{align*}
C_{\Gamma}^{3,1}=&\mathbb{C}\{dz^{123\bar{1}},e^{-z^{1}}dz^{123\bar{2}},e^{z^{1}}dz^{123\bar{3}},e^{-z^{\bar{1}}}dz^{123\bar{2}},e^{z^{\bar{1}}}dz^{123\bar{3}} \}\\
C_{\Gamma}^{2,2}=&\mathbb{C}\{e^{-2z^{1}}dz^{12\overline{12}},dz^{12\overline{13}},e^{-z^{1}}dz^{12\overline{23}},dz^{13\overline{12}},e^{2z^{1}}dz^{13\overline{13}},e^{z^{1}}dz^{13\overline{23}},
e^{-z^{1}}dz^{23\overline{12}},\\
&e^{z^{1}}dz^{23\overline{13}},dz^{23\overline{23}},e^{-2z^{\bar{1}}}dz^{12\overline{12}},e^{-z^{\bar{1}}}dz^{23\overline{12}},e^{-z^{\bar{1}}}dz^{12\overline{23}},e^{z^{\bar{1}}}dz^{13\overline{23}},e^{2z^{\bar{1}}}dz^{13\overline{13}}, e^{z^{\bar{1}}}dz^{23\overline{13}}\}\\
C_{\Gamma}^{1,3}=&\mathbb{C}\{dz^{1\overline{123}},e^{-z^{\bar{1}}}dz^{2\overline{123}},e^{z^{\bar{1}}}dz^{3\overline{123}},e^{-z^{1}}dz^{2\overline{123}},e^{z^{1}}dz^{3\overline{123}}\}\\
C_{\Gamma}^{3,2}=&\mathbb{C}\{e^{-z^{1}}dz^{123\overline{12}},e^{z^{1}}dz^{123\overline{13}},dz^{123\overline{23}},e^{-z^{\bar{1}}}dz^{123\overline{12}},e^{z^{\bar{1}}}dz^{123\overline{13}}\}\\
C_{\Gamma}^{2,3}=&\mathbb{C}\{e^{-z^{1}}dz^{12\overline{123}},e^{z^{1}}dz^{13\overline{123}},dz^{23\overline{123}},e^{-z^{\bar{1}}}dz^{12\overline{123}},e^{z^{\bar{1}}}dz^{13\overline{123}}\}
\end{align*}
By using \eqref{eq-Cgamma-A}, we easily get the following Aeppli cohomology:
\begin{align*}
H_{A}^{1,0}(X)=&\mathbb{C}\{dz^1,e^{-z^1}dz^2,e^{z^1}dz^3,e^{-z^{\bar{1}}}dz^2,e^{z^{\bar{1}}}dz^3\}\\
H_{A}^{0,1}(X)=&\mathbb{C}\{dz^{\bar{1}},e^{-z^1}dz^{\bar{2}},e^{z^1}dz^{\bar{3}},e^{-z^{\bar{1}}}dz^{\bar{2}},e^{z^{\bar{1}}}dz^{\bar{3}}\}\\
H_{A}^{1,1}(X)=&\mathbb{C}\{dz^{1\bar{1}},e^{-z^{1}}dz^{2\bar{1}},e^{-2z^{1}}dz^{2\bar{2}},dz^{2\bar{3}},e^{z^{1}}dz^{3\bar{1}},dz^{3\bar{2}},\\
&e^{2z^{1}}dz^{3\bar{3}},e^{-z^{\bar{1}}}dz^{1\bar{2}},e^{z^{\bar{1}}}dz^{1\bar{3}},e^{-2z^{\bar{1}}}dz^{2\bar{2}},e^{2z^{\bar{1}}}dz^{3\bar{3}}, \}\\
H_{A}^{2,1}(X)=&\mathbb{C}\{dz^{12\bar{3}},dz^{13\bar{2}},dz^{23\bar{1}},e^{-z^{1}}dz^{23\bar{2}},e^{z^{1}}dz^{23\bar{3}},e^{-2z^{\bar{1}}}dz^{12\bar{2}},e^{2z^{\bar{1}}}dz^{13\bar{3}},e^{-z^{\bar{1}}}dz^{23\bar{2}},e^{z^{\bar{1}}}dz^{23\bar{3}} \}\\
H_{A}^{1,2}(X)=&\mathbb{C}\{dz^{1\overline{23}}, e^{-2z^{1}}dz^{2\overline{12}},dz^{2\overline{13}},e^{-z^{1}}dz^{2\overline{23}},dz^{3\overline{12}},e^{2z^{1}}dz^{3\overline{13}},e^{z^{1}}dz^{3\overline{23}},e^{-z^{\bar{1}}}dz^{2\overline{23}},e^{z^{\bar{1}}}dz^{3\overline{23}} \}
\end{align*}
$X$ is equipped with the Hermitian metric $h=\sum_{i=1}^3\varphi^i\otimes\bar{\varphi}^i$ where
\[
\varphi^1=dz^1,\quad\varphi^2=e^{-z^1}dz^2,\quad\varphi^3=e^{z^1}dz^3.
\]
Under this metric, the representatives of these cohomology classes listed as above are all harmonic.

Let us first consider canonical Aeppli deformations.

For $(p,q)=(1,0)$. Let
\[
\sigma_0 = a_{1}dz^1+a_{2}e^{-z^1}dz^2+a_{3}e^{z^1}dz^3+a_{4}e^{-z^{\bar{1}}}dz^2+a_{5}e^{z^{\bar{1}}}dz^3\in\mathcal{H}_{A}^{1,0}(X),
\]
we have
\begin{align*}
\p\mathcal{L}_{\phi_1}^{1,0}\sigma_0&=\p[-a_1\p i_{\phi_{1}}dz^{1}+a_{2}i_{\phi_{1}}\p(e^{-z^{1}}dz^{2})+a_{3}i_{\phi_{1}}\p(e^{z^{1}}dz^{3})]\\
&=\p[-a_{1}\p(tdz^{\bar{1}})-a_{2}e^{-z^1}i_{\phi_{1}}dz^{12}+a_{3}e^{z^1}i_{\phi_{1}}dz^{13}]\\
&=\p[a_2te^{-z^{1}}dz^{2\bar{1}}-a_{3}te^{z^{1}}dz^{3\bar{1}}]\\
&=-a_{2}te^{-z^{1}}dz^{12\bar{1}}-a_{3}te^{z^{1}}dz^{13\bar{1}},
\end{align*}
which is $\p\pb$-exact if and only if $\p\mathcal{L}_{\phi_1}^{1,0}\sigma_0=0$, i.e.
\begin{equation}\label{eq-A1,0}
a_{2}t=a_{3}t=0
\end{equation}
and in this case the canonical Aeppli deformation is given by
\[
\sigma_1 = G_{A}(\p\bar{\p})^* \p i_{\phi_1}\p \sigma_0=0,\quad  \sigma_2 = G_{A}(\p\bar{\p})^* (\p i_{\phi_2}\p \sigma_0 + \p i_{\phi_1}\p \sigma_1 )=0.
\]
Similarly, there holds
\[
\p\sum_{j=1}^{k}\mathcal{L}_{\phi_j}^{1,0}\sigma_{k-j}=0 , \quad k\geq 2,
\]
and $\sigma_k= G_{A}(\p\bar{\p})^*\sum_{i+j=k} \p i_{\phi_j} \p\sigma_i=0$ for any $k\geq2$. Therefore, for $W=\mathcal{H}_{A}^{1,0}(X)$ we have (See Definition \ref{def-W_tp,q-f_t})
\begin{align*}
W_{t}^{1,0}=
&\Big\{ \sigma_0\in W \mid  \sigma(t)\in\ker \p\pb_{\phi(t)},\\
&\text{where}~\sigma (t)=\sum_{k} \sigma_k~\text{with}~\sigma_k=G_{A}(\p\bar{\p})^*\sum_{i+j=k} \p i_{\phi_j} \p\sigma_i,~\forall k>0 \Big\}\\
=&\Big\{ \sigma_0\in W \mid  \eqref{eq-A1,0}~\text{is~satisfied} \Big\},
\end{align*}
and
\[
w_t^{1,0}=\dim H_{A}^{1,0}(X) - \dim W_{t}^{1,0}=
\left\{
\begin{array}{rcl}
0, &t =0 \\[5pt]
2, &t\neq 0.
\end{array}
\right.
\]

For $(p,q)=(0,1)$. Let
\[
\sigma_0 = a_1dz^{\bar{1}}+a_2e^{-z^1}dz^{\bar{2}}+a_3e^{z^1}dz^{\bar{3}}+a_4e^{-z^{\bar{1}}}dz^{\bar{2}}+a_5e^{z^{\bar{1}}}dz^{\bar{3}}\in\mathcal{H}_{A}^{0,1}(X),
\]
we have
\begin{align*}
 \p\mathcal{L}_{\phi_1}^{1,0}\sigma_0&=\p[a_2 i_{\phi_{1}}\p(e^{-z^{1}}dz^{\bar{2}})+a_3 i_{\phi_{1}}\p(e^{z^{1}}dz^{\bar{3}})]\\
 &=\p(-a_2te^{-z^{1}}dz^{\overline{12}}+a_3te^{z^{1}}dz^{\overline{13}})\\
 &=a_2te^{-z^{1}}dz^{1\overline{12}}+a_3te^{z^{1}}dz^{1\overline{13}}
\end{align*}
is $\p\pb$-exact if and only if $\p\mathcal{L}_{\phi_1}^{1,0}\sigma_0=0$, i.e.
\begin{equation}\label{eq-A0,1}
a_2t=a_3t=0
\end{equation}
and in this case $\sigma_k= G_{A}(\p\bar{\p})^*\sum_{i+j=k} \p i_{\phi_j} \p\sigma_i=0$ for any $k\geq1$. Therefore, for $W=\mathcal{H}_{A}^{0,1}(X)$ we have
\begin{align*}
W_{t}^{0,1}=\Big\{ \sigma_0\in W \mid  \eqref{eq-A0,1}~\text{is~satisfied} \Big\},
\end{align*}
and
\[
w_t^{0,1}=\dim H_{A}^{0,1}(X) - \dim W_{t}^{0,1}=
\left\{
\begin{array}{rcl}
0, &t =0 \\[5pt]
2, &t\neq 0.
\end{array}
\right.
\]

For $(p,q)=(1,1)$. Let
\begin{align*}
\sigma_0 =& a_1dz^{1\bar{1}}+a_2e^{-z^{1}}dz^{2\bar{1}}+a_3e^{-2z^{1}}dz^{2\bar{2}}+a_4dz^{2\bar{3}}+a_5e^{z^{1}}dz^{3\bar{1}}+\\
&a_6dz^{3\bar{2}}+a_7e^{2z^{1}}dz^{3\bar{3}}+a_8e^{-z^{\bar{1}}}dz^{1\bar{2}}+a_9e^{z^{\bar{1}}}dz^{1\bar{3}}+a_{10}e^{-2z^{\bar{1}}}dz^{2\bar{2}}+a_{11}e^{2z^{\bar{1}}}dz^{3\bar{3}}\in\mathcal{H}_{A}^{1,1}(X),
\end{align*}
we have
\begin{align*}
\p\mathcal{L}_{\phi_1}^{1,0}\sigma_0&=\p[a_3 i_{\phi_{1}}\p(e^{-2z^{1}}dz^{2\bar{2}})+a_7 i_{\phi_{1}}\p(e^{2z^{1}}dz^{3\bar{3}})]\\
&=\p(2a_{3}te^{-2z^{1}}dz^{2\overline{12}}-2a_{7}te^{2z^{1}}dz^{3\overline{13}})\\
&=-4a_{3}te^{-2z^{1}}dz^{12\overline{12}}-4a_{7}te^{2z^{1}}dz^{13\overline{13}}
\end{align*}
is $\p\pb$-exact if and only if $\p\mathcal{L}_{\phi_1}^{1,0}\sigma_0=0$, i.e.
\begin{equation}\label{eq-A1,1}
a_3t=a_7t=0
\end{equation}
and in this case $\sigma_k= G_{A}(\p\bar{\p})^*\sum_{i+j=k} \p i_{\phi_j} \p\sigma_i=0$ for any $k\geq1$. Therefore, for $W=\mathcal{H}_{A}^{1,1}(X)$ we have
\begin{align*}
W_{t}^{1,1}=\Big\{ \sigma_0\in W \mid  \eqref{eq-A1,1}~\text{is~satisfied} \Big\},
\end{align*}
and
\[
w_t^{1,1}=\dim H_{A}^{1,1}(X) - \dim W_{t}^{1,1}=
\left\{
\begin{array}{rcl}
0, &t =0 \\[5pt]
2, &t\neq 0.
\end{array}
\right.
\]

For $(p,q)=(2,1)$. Let
\begin{align*}
\sigma_0 =& a_1dz^{12\bar{3}}+a_2dz^{13\bar{2}}+a_3dz^{23\bar{1}}+a_4e^{-z^{1}}dz^{23\bar{2}}+\\
&a_5e^{z^{1}}dz^{23\bar{3}}+a_6e^{-2z^{\bar{1}}}dz^{12\bar{2}}+a_7e^{2z^{\bar{1}}}dz^{13\bar{3}}+a_8e^{-z^{\bar{1}}}dz^{23\bar{2}}+a_9e^{z^{\bar{1}}}dz^{23\bar{3}}\in\mathcal{H}_{A}^{2,1}(X),
\end{align*}
we have
\begin{align*}
\p\mathcal{L}_{\phi_1}^{1,0}\sigma_0&=\p[a_4i_{\phi_{1}}\p(e^{-z^{1}}dz^{23\bar{2}})+a_5i_{\phi_{1}}\p(e^{z^{1}}dz^{23\bar{3}})\\
&-a_6\p i_{\phi_{1}}(e^{-2z^{\bar{1}}}dz^{12\bar{2}})-a_7\p i_{\phi_{1}}(e^{2z^{\bar{1}}}dz^{13\bar{3}})]\\
&=\p(-a_4te^{-z^{1}}dz^{23\overline{12}}+a_5te^{z^{1}}dz^{23\overline{13}})\\
&=a_4te^{-z^{1}}dz^{123\overline{12}}+a_5te^{z^{1}}dz^{123\overline{13}}
\end{align*}
is $\p\pb$-exact if and only if $\p\mathcal{L}_{\phi_1}^{1,0}\sigma_0=0$, i.e.
\begin{equation}\label{eq-A2,1}
a_{4}t=a_{5}t=0
\end{equation}
and in this case $\sigma_k= G_{A}(\p\bar{\p})^*\sum_{i+j=k} \p i_{\phi_j} \p\sigma_i=0$ for any $k\geq1$. Therefore, for $W=\mathcal{H}_{A}^{2,1}(X)$ we have
\begin{align*}
W_{t}^{2,1}=\Big\{ \sigma_0\in W \mid  \eqref{eq-A2,1}~\text{is~satisfied} \Big\},
\end{align*}
and
\[
w_t^{2,1}=\dim H_{A}^{2,1}(X) - \dim W_{t}^{2,1}=
\left\{
\begin{array}{rcl}
0, &t =0 \\[5pt]
2, &t\neq 0.
\end{array}
\right.
\]

For $(p,q)=(1,2)$. Let
\begin{align*}
\sigma_0 = &a_1dz^{1\overline{23}}+a_2 e^{-2z^{1}}dz^{2\overline{12}}+a_3dz^{2\overline{13}}+a_4e^{-z^{1}}dz^{2\overline{23}}+\\
&a_5dz^{3\overline{12}}+a_6e^{2z^{1}}dz^{3\overline{13}}+a_7e^{z^{1}}dz^{3\overline{23}}+a_8e^{-z^{\bar{1}}}dz^{2\overline{23}}+a_9e^{z^{\bar{1}}}dz^{3\overline{23}}\in\mathcal{H}_{A}^{1,2}(X),
\end{align*}
we have
\begin{align*}
\p\mathcal{L}_{\phi_1}^{1,0}\sigma_0&=\p[a_4i_{\phi_{1}}\p(e^{-z^{1}}dz^{2\overline{23}})+a_7i_{\phi_{1}}\p(e^{z^{1}}dz^{3\overline{23}})]\\
&=\p(a_4te^{-z^{1}}dz^{2\overline{123}}-a_7te^{z^{1}}dz^{3\overline{123}})\\
&=-a_4te^{-z^{1}}dz^{12\overline{123}}-a_7te^{z^{1}}dz^{13\overline{123}}
\end{align*}
is $\p\pb$-exact if and only if $\p\mathcal{L}_{\phi_1}^{1,0}\sigma_0=0$, i.e.
\begin{equation}\label{eq-A1,2}
a_{4}t=a_{7}t=0
\end{equation}
and in this case $\sigma_k= G_{A}(\p\bar{\p})^*\sum_{i+j=k} \p i_{\phi_j} \p\sigma_i=0$ for any $k\geq1$. Therefore, for $W=\mathcal{H}_{A}^{1,2}(X)$ we have
\begin{align*}
W_{t}^{1,2}=\Big\{ \sigma_0\in W \mid  \eqref{eq-A1,2}~\text{is~satisfied} \Big\},
\end{align*}
and
\[
w_t^{1,2}=\dim H_{A}^{1,2}(X) - \dim W_{t}^{1,2}=
\left\{
\begin{array}{rcl}
0, &t =0 \\[5pt]
2, &t\neq 0.
\end{array}
\right.
\]

\vspace{1cm}
Now we consider canonical Bott-Chern deformations.

By using \eqref{eq-Cgamma-BC}, we get
\begin{align*}
H_{BC}^{2,0}(X)=&\mathbb{C}\{e^{-z^{1}}dz^{12},e^{z^{1}}dz^{13},dz^{23} \}\\
H_{BC}^{1,1}(X)=&\mathbb{C}\{dz^{1\bar{1}},e^{-z^{1}}dz^{1\bar{2}},e^{z^{1}}dz^{1\bar{3}},dz^{2\bar{3}},dz^{3\bar{2}},e^{-z^{\bar{1}}}dz^{2\bar{1}},e^{z^{\bar{1}}}dz^{3\bar{1}} \}\\
H_{BC}^{2,1}(X)=&\mathbb{C}\{e^{-z^{1}}dz^{12\bar{1}},e^{-2z^{1}}dz^{12\bar{2}},dz^{12\bar{3}},e^{z^{1}}dz^{13\bar{1}},dz^{13\bar{2}},e^{2z^{1}}dz^{13\bar{3}},dz^{23\bar{1}},e^{-z^{\bar{1}}}dz^{12\bar{1}},e^{z^{\bar{1}}}dz^{13\bar{1}}  \}\\
H_{BC}^{1,2}(X)=&\mathbb{C}\{ e^{-z^{1}}dz^{1\overline{12}},e^{z^{1}}dz^{1\overline{13}},dz^{1\overline{23}},dz^{2\overline{13}},dz^{3\overline{12}},e^-{z^{\bar{1}}}dz^{1\overline{12}},
e^{z^{\bar{1}}}dz^{1\overline{13}},e^{-2z^{\bar{1}}}dz^{2\overline{12}},e^{2z^{\bar{1}}}dz^{3\overline{13}} \}\\
H_{BC}^{3,1}(X)=&\mathbb{C}\{dz^{123\bar{1}},e^{-z^{1}}dz^{123\bar{2}},e^{z^{1}}dz^{123\bar{3}} \}\\
H_{BC}^{2,2}(X)=&\mathbb{C}\{e^{-2z^{1}}dz^{12\overline{12}},dz^{12\overline{13}},e^{-z^{1}}dz^{12\overline{23}},dz^{13\overline{12}},e^{2z^{1}}dz^{13\overline{13}},e^{z^{1}}dz^{13\overline{23}},\\
&dz^{23\overline{23}},e^{-2z^{\bar{1}}}dz^{12\overline{12}},e^{2z^{\bar{1}}}dz^{13\overline{13}},e^{-z^{\bar{1}}}dz^{23\overline{12}}, e^{z^{\bar{1}}}dz^{23\overline{13}}\}\\
H_{BC}^{3,2}(X)=&\mathbb{C}\{e^{-z^{1}}dz^{123\overline{12}},e^{z^{1}}dz^{123\overline{13}},dz^{123\overline{23}},e^{-z^{\bar{1}}}dz^{123\overline{12}},e^{z^{\bar{1}}}dz^{123\overline{13}}\}\\
H_{BC}^{2,3}(X)=&\mathbb{C}\{e^{-z^{1}}dz^{12\overline{123}},e^{z^{1}}dz^{13\overline{123}},dz^{23\overline{123}},e^{-z^{\bar{1}}}dz^{12\overline{123}},e^{z^{\bar{1}}}dz^{13\overline{123}}\}
\end{align*}
For $(p,q)=(2,0)$. Let $\sigma_0 = a_1e^{-z^{1}}dz^{12}+a_2e^{z^{1}}dz^{13}+a_3dz^{23}\in\mathcal{H}_{BC}^{2,0}(X)$, we have
\begin{align*}
 \p i_{\phi_1}\sigma_0&=-a_1t\p(e^{-z^{1}}dz^{2\bar{1}})-a_2t\p(e^{z^{1}}dz^{3\bar{1}})\\
 &=a_1te^{-z^{1}}dz^{12\bar{1}}-a_2te^{z^{1}}dz^{13\bar{1}}
\end{align*}
is $\pb$-exact if and only if $ \p i_{\phi_1}\sigma_0=0$, i.e.
\begin{equation}\label{eq-BC2,0}
a_1t=a_2t=0
\end{equation}
and in this case
\[
\sigma_{1}=-G_{BC}(\pb^*\p\p^*+\pb^*)\partial i_{\phi_1}\sigma_0 =0.
\]
On the other hand,
\[
\partial i_{\phi_2}\sigma_0 = 0 \Longrightarrow \sigma_{2}=-G_{BC}(\pb^*\p\p^*+\pb^*)\partial(i_{\phi_2}\sigma_0+i_{\phi_1}\sigma_1) =0,
\]
and $\phi_k=0,~ k\geq2$  implies that $\sigma_{k}= G_{BC}(\pb^*\p\p^*+\pb^*)\sum_{i+j=k} \p i_{\phi_j}\sigma_i=0$ for any $k\geq2$.

Therefore, for $V=\mathcal{H}_{BC}^{2,0}(X)$ we have (See \eqref{eq-V BC,t})
\[
V_{BC,t}^{2,0}=\{a_1e^{-z^{1}}dz^{12}+a_2e^{z^{1}}dz^{13}+a_3dz^{23} \mid (a_{1},a_2.a_3)\in \mathbb{C}^3~\text{satisfy}~\eqref{eq-BC2,0}\},
\]
and
\[
v_t^{2,0}=\dim H_{BC}^{2,0}(X)-\dim V_{BC,t}^{2,0}=
\left\{
\begin{array}{rcl}
0, &t =0 \\[5pt]
2, &t \neq 0.
\end{array}
\right.
\]

For $(p,q)=(1,1)$. Let
\[
\sigma_0 = a_1dz^{1\bar{1}}+a_2e^{-z^{1}}dz^{1\bar{2}}+a_3e^{z^{1}}dz^{1\bar{3}}+a_4dz^{2\bar{3}}+
a_5dz^{3\bar{2}}+a_6e^{-z^{\bar{1}}}dz^{2\bar{1}}+a_7e^{z^{\bar{1}}}dz^{3\bar{1}}\in\mathcal{H}_{BC}^{1,1}(X),
\]
we have
\begin{align*}
 \p i_{\phi_1}\sigma_0&=a_2t\p(e^{-z^{1}}dz^{\overline{12}})+a_3t\p(e^{z^{1}}dz^{\overline{13}})\\
 &=-a_2te^{-z^{1}}dz^{1\overline{12}}+a_3te^{z^{1}}dz^{1\overline{13}}
\end{align*}
is $\pb$-exact if and only if $ \p i_{\phi_1}\sigma_0=0$, i.e.
\begin{equation}\label{eq-BC1,1}
a_2t=a_3t=0
\end{equation}
and in this case $\sigma_{k}= G_{BC}(\pb^*\p\p^*+\pb^*)\sum_{i+j=k} \p i_{\phi_j}\sigma_i=0$ for any $k\geq1$.

Therefore, for $V=\mathcal{H}_{BC}^{1,1}(X)$ we have
\begin{align*}
 V_{BC,t}^{1,1}&=\{a_1dz^{1\bar{1}}+a_2e^{z^{1}}dz^{1\bar{2}}+a_3e^{-z^{1}}dz^{1\bar{3}}+a_4dz^{2\bar{3}}+a_5dz^{3\bar{2}}+a_6e^{z^{\bar{1}}}dz^{2\bar{1}}\\
 &+a_7e^{-z^{\bar{1}}}dz^{3\bar{1}} \mid (a_{1},\cdots,a_7)\in \mathbb{C}^7~\text{satisfy}~\eqref{eq-BC1,1}\},
\end{align*}
and
\[
v_t^{1,1}=\dim H_{BC}^{1,1}(X)-\dim V_{BC,t}^{1,1}=
\left\{
\begin{array}{rcl}
0, &t =0 \\[5pt]
2, &t \neq 0.
\end{array}
\right.
\]

For $(p,q)=(2,1)$. Let
\begin{align*}
\sigma_0 =& a_1e^{-z^{1}}dz^{12\bar{1}}+a_2e^{-2z^{1}}dz^{12\bar{2}}+a_3dz^{12\bar{3}}+a_4e^{z^{1}}dz^{13\bar{1}}+a_5dz^{13\bar{2}}+\\
&a_6e^{2z^{1}}dz^{13\bar{3}}+a_7dz^{23\bar{1}}+a_8e^{-z^{\bar{1}}}dz^{12\bar{1}}+a_9e^{z^{\bar{1}}}dz^{13\bar{1}}\in\mathcal{H}_{BC}^{2,1}(X),
\end{align*}
we have
\begin{align*}
 \p i_{\phi_1}\sigma_0&=-a_2t\p(e^{-2z^{1}}dz^{2\overline{12}})-a_6t\p(e^{2z^{1}}dz^{3\overline{13}})\\
&=2a_2te^{-2z^{1}}dz^{12\overline{12}}-2a_6te^{2z^{1}}dz^{13\overline{13}}
\end{align*}
is $\pb$-exact if and only if $\p i_{\phi_1}\sigma_0=0$, i.e.
\begin{equation}\label{eq-BC2,1}
a_2t=a_6t=0
\end{equation}
and in this case $\sigma_{k}= G_{BC}(\pb^*\p\p^*+\pb^*)\sum_{i+j=k} \p i_{\phi_j}\sigma_i=0$ for any $k\geq1$.

Therefore, for $V=\mathcal{H}_{BC}^{2,1}(X)$ we have
\begin{align*}
V_{BC,t}^{2,1}&=\{a_1e^{z^{1}}dz^{12\bar{1}}+a_2e^{2z^{1}}dz^{12\bar{2}}+a_3dz^{12\bar{3}}+a_4e^{-z^{1}}dz^{13\bar{1}}+a_5dz^{13\bar{2}}+a_6e^{-2z^{1}}dz^{13\bar{3}}\\
&+a_7dz^{23\bar{1}}+a_8e^{z^{\bar{1}}}dz^{12\bar{1}}+a_9e^{-z^{\bar{1}}}dz^{13\bar{1}} \mid (a_{1},\cdots,a_9)\in \mathbb{C}^9~\text{satisfy}~\eqref{eq-BC2,1}\},
\end{align*}
and
\[
v_t^{2,1}=\dim H_{BC}^{2,1}(X)-\dim V_{BC,t}^{2,1}=
\left\{
\begin{array}{rcl}
0, &t =0 \\[5pt]
2, &t \neq 0.
\end{array}
\right.
\]

For $(p,q)=(3,1)$. Let
\[
\sigma_0 = a_1dz^{123\bar{1}}+a_2e^{-z^{1}}dz^{123\bar{2}}+a_3e^{z^{1}}dz^{123\bar{3}}\in\mathcal{H}_{BC}^{3,1}(X),
\]
we have
\begin{align*}
 \p i_{\phi_1}\sigma_0&=a_2t\p(-e^{z^{1}}dz^{23\overline{12}})+a_3t\p(e^{z^{1}}dz^{23\overline{13}})\\
&=-a_2te^{-z^{1}}dz^{123\overline{12}}+a_3te^{z^{1}}dz^{123\overline{13}}
\end{align*}
is $\pb$-exact if and only if $\p i_{\phi_1}\sigma_0=0$, i.e.
\begin{equation}\label{eq-BC3,1}
a_2t=a_3t=0
\end{equation}
and in this case $\sigma_{k}= G_{BC}(\pb^*\p\p^*+\pb^*)\sum_{i+j=k} \p i_{\phi_j}\sigma_i=0$ for any $k\geq1$.

Therefore, for $V=\mathcal{H}_{BC}^{3,1}(X)$ we have
\begin{align*}
 V_{BC,t}^{3,1}=\{a_1dz^{123\bar{1}}+a_2e^{z^{1}}dz^{123\bar{2}}+a_3e^{-z^{1}}dz^{123\bar{3}} \mid (a_{1},a_2,a_3)\in \mathbb{C}^3~\text{satisfy}~\eqref{eq-BC3,1}\},
\end{align*}
and
\[
v_t^{3,1}=\dim H_{BC}^{3,1}(X)-\dim V_{BC,t}^{3,1}=
\left\{
\begin{array}{rcl}
0, &t =0 \\[5pt]
2, &t \neq 0.
\end{array}
\right.
\]

For $(p,q)=(2,2)$. Let
\begin{align*}
\sigma_0 =& a_1e^{-2z^{1}}dz^{12\overline{12}}+a_2dz^{12\overline{13}}+a_3e^{-z^{1}}dz^{12\overline{23}}+a_4dz^{13\overline{12}}+a_5e^{2z^{1}}dz^{13\overline{13}}+a_6e^{z^{1}}dz^{13\overline{23}}+\\
&a_7dz^{23\overline{23}}+a_8e^{-2z^{\bar{1}}}dz^{12\overline{12}}+a_9e^{2z^{\bar{1}}}dz^{13\overline{13}}+a_{10}e^{-z^{\bar{1}}}dz^{23\overline{12}}+a_{11}e^{z^{\bar{1}}}dz^{23\overline{13}}\in\mathcal{H}_{BC}^{2,2}(X),
\end{align*}
we have
\begin{align*}
 \p i_{\phi_1}\sigma_0&=-a_3t\p(e^{-z^{1}}dz^{2\overline{123}})-a_6t\p(e^{z^{1}}dz^{3\overline{123}})\\
&=a_3te^{-z^{1}}dz^{12\overline{123}}-a_6te^{z^{1}}dz^{13\overline{123}}
\end{align*}
is $\pb$-exact if and only if $\p i_{\phi_1}\sigma_0=0$, i.e.
\begin{equation}\label{eq-BC2,2}
a_3t=a_6t=0
\end{equation}
and in this case $\sigma_{k}= G_{BC}(\pb^*\p\p^*+\pb^*)\sum_{i+j=k} \p i_{\phi_j}\sigma_i=0$ for any $k\geq1$.

Therefore, for $V=\mathcal{H}_{BC}^{2,2}(X)$ we have
\begin{align*}
 V_{BC,t}^{2,2}&=\{a_1e^{2z^{1}}dz^{12\overline{12}}+a_2dz^{12\overline{13}}+a_3e^{z^{1}}dz^{12\overline{23}}+a_4dz^{13\overline{12}}+a_5e^{-2z^{1}}dz^{13\overline{13}}\\
 &+a_6e^{-z^{1}}dz^{13\overline{23}}+a_7dz^{23\overline{23}}+a_8e^{2z^{\bar{1}}}dz^{12\overline{12}}+a_9e^{-2z^{\bar{1}}}dz^{13\overline{13}}+a_{10}e^{z^{\bar{1}}}dz^{23\overline{12}}\\
 &+a_{11}e^{-z^{\bar{1}}}dz^{23\overline{13}} \mid (a_{1},\cdots,a_{11})\in \mathbb{C}^{11}~\text{satisfy}~\eqref{eq-BC2,2}\},
\end{align*}
and
\[
v_t^{2,2}=\dim H_{BC}^{2,2}(X)-\dim V_{BC,t}^{2,2}=
\left\{
\begin{array}{rcl}
0, &t =0 \\[5pt]
2, &t \neq 0.
\end{array}
\right.
\]

Similarly, we can determine $w_t^{p,q}$ and $v_t^{p,q}$ for all $(p,q)$. Collecting computations in this subsection all together, we get the following table (in the next page):
\renewcommand\arraystretch{1.5}
\begin{table}[!htbp]
\caption{$v^{\bullet,\bullet}_t$ and $w^{\bullet,\bullet}_t$ for $t\neq0$}
\centering
\begin{center}
\begin{tabular}{|c|c|c|c|c|c|c|c|c|c|c|c|}
\hline
    $h^{0,1}_A$ & $5$ & $h^{3,2}_{BC}$ & $5$ & $v^{3,2}_t$ & $0$ & $w^{2,1}_t$ & $2$ & $h^{3,2}_{BC,\phi(t)}$ & $3$ & $h^{0,1}_{A,\phi(t)}$ & $3$  \\
\hline
   $h^{1,0}_A$ & $5$ & $h^{2,3}_{BC}$ & $5$ & $v^{2,3}_t$ & $0$ & $w^{1,2}_t$ & $2$ & $h^{2,3}_{BC,\phi(t)}$ & $3$ & $h^{1,0}_{A,\phi(t)}$ & $3$ \\
\hline
   $h^{0,2}_A$ & $3$ & $h^{3,1}_{BC}$ & $3$ & $v^{3,1}_t$ & $2$ & $w^{2,0}_t$ & $0$ & $h^{3,1}_{BC,\phi(t)}$ & $1$ & $h^{0,2}_{A,\phi(t)}$ & $1$  \\
\hline
   $h^{1,1}_A$ & $11$ & $h^{2,2}_{BC}$ & $11$ & $v^{2,2}_t$ & $2$ & $w^{1,1}_t$ & $2$ & $h^{2,2}_{BC,\phi(t)}$ & $7$ & $h^{1,1}_{A,\phi(t)}$ & $7$  \\
\hline
  $h^{2,0}_A$ & $3$ & $h^{1,3}_{BC}$ & $3$ & $v^{1,3}_t$ & $0$ & $w^{0,2}_t$ & $0$ & $h^{1,3}_{BC,\phi(t)}$ & $3$ & $h^{2,0}_{A,\phi(t)}$ & $3$  \\
\hline
  $h^{0,3}_A$ & $1$ & $h^{3,0}_{BC}$ & $1$ & $v^{3,0}_t$ & $0$ & $w^{2,-1}_t$ & $0$ & $h^{3,0}_{BC,\phi(t)}$ & $1$ & $h^{0,3}_{A,\phi(t)}$ & $1$  \\
\hline
   $h^{1,2}_A$ & $9$ & $h^{2,1}_{BC}$ & $9$ & $v^{2,1}_t$ & $2$ & $w^{1,0}_t$ & $2$ & $h^{2,1}_{BC,\phi(t)}$ & $5$ & $h^{1,2}_{A,\phi(t)}$ & $5$  \\
 \hline
    $h^{2,1}_A$ & $9$ & $h^{1,2}_{BC}$ & $9$ & $v^{1,2}_t$ & $0$ & $w^{0,1}_t$ & $2$ & $h^{1,2}_{BC,\phi(t)}$ & $7$ & $h^{2,1}_{A,\phi(t)}$ & $7$ \\
\hline
   $h^{3,0}_A$ & $1$ & $h^{0,3}_{BC}$ & $1$ & $v^{0,3}_t$ & $0$ & $w^{-1,2}_t$ & $0$ & $h^{0,3}_{BC,\phi(t)}$ & $1$ & $h^{3,0}_{A,\phi(t)}$ & $1$  \\
\hline
   $h^{1,3}_A$ & $3$ & $h^{2,0}_{BC}$ & $3$ & $v^{2,0}_t$ & $2$ & $w^{1,-1}_t$ & $0$ & $h^{2,0}_{BC,\phi(t)}$ & $1$ & $h^{1,3}_{A,\phi(t)}$ & $1$ \\
\hline
   $h^{2,2}_A$ & $7$ & $h^{1,1}_{BC}$ & $7$ & $v^{1,1}_t$ & $2$ & $w^{0,0}_t$ & $0$ & $h^{1,1}_{BC,\phi(t)}$ & $5$ & $h^{2,2}_{A,\phi(t)}$ & $5$  \\
\hline
   $h^{3,1}_A$ & $3$ & $h^{0,2}_{BC}$ & $3$ & $v^{0,2}_t$ & $0$ & $w^{-1,1}_t$ & $0$ & $h^{0,2}_{BC,\phi(t)}$ & $3$ & $h^{3,1}_{A,\phi(t)}$ & $3$  \\
\hline
   $h^{2,3}_A$ & $1$ & $h^{1,0}_{BC}$ & $1$ & $v^{1,0}_t$ & $0$ & $w^{0,-1}_t$ & $0$ & $h^{1,0}_{BC,\phi(t)}$ & $1$ & $h^{2,3}_{A,\phi(t)}$ & $1$  \\
\hline
    $h^{3,2}_A$ & $1$ & $h^{0,1}_{BC}$ & $1$ & $v^{0,1}_t$ & $0$ & $w^{-1,0}_t$ & $0$ & $h^{0,1}_{BC,\phi(t)}$ & $1$ & $h^{3,2}_{A,\phi(t)}$ & $1$  \\
\hline
\end{tabular}
\end{center}
\end{table}

\bibliographystyle{alpha}
\bibliography{reference}
\end{document}